\documentclass{amsart}
\usepackage{graphicx}
\usepackage{amsmath}
\usepackage{amscd}
\usepackage{amsfonts}
\usepackage{amssymb}
\newtheorem{theorem}{Theorem}[section]
\theoremstyle{plain}

\newtheorem{corollary}[theorem]{Corollary}

\newtheorem{definition}[theorem]{Definition}

\newtheorem{lemma}[theorem]{Lemma}

\newtheorem{problem}[theorem]{Problem}
\newtheorem{proposition}[theorem]{Proposition}
\newtheorem{remark}[theorem]{Remark}

\numberwithin{equation}{section}

\begin{document}
\title[Baire-1 Functions]{Functions of Baire Class One}
\author{Denny H. Leung}
\address{Department of Mathematics, National University of Singapore, 2 Science Drive
2, Singapore 117543.}
\email{matlhh@nus.edu.sg}
\author{Wee-Kee Tang}
\address{Division of Mathematics, National Institute of Education \\
Nanyang Technological University, 469 Bukit Timah Road, Singapore 259756.}
\email{wktang@nie.edu.sg}
\begin{abstract}Let $K$ be a compact metric space. A real-valued function on $K$ is said to be
of Baire class one (Baire-$1$) if it is the pointwise limit of a sequence of
continuous functions. In this paper, we study two well known ordinal indices
of Baire-$1$ functions, the oscillation index $\beta$ and the convergence
index $\gamma$. It is shown that these two indices are fully compatible in the
following sense : a Baire-$1$ function $f$ satisfies $\beta(f)\leq\omega
^{\xi_{1}}\cdot\omega^{\xi_{2}}$ for some countable ordinals $\xi_{1}$ and
$\xi_{2}$ if and only if there exists a sequence of Baire-$1$ functions
$(f_{n})$ converging to $f$ pointwise such that $\sup_{n}\beta(f_{n}%
)\leq\omega^{\xi_{1}}$ and $\gamma((f_{n}))\leq\omega^{\xi_{2}}$. We also
obtain an extension result for Baire-$1$ functions analogous to the Tietze
Extension Theorem. Finally, it is shown that if $\beta\left(  f\right)
\leq\omega^{\xi_{1}}$ and $\beta\left(  g\right)  \leq\omega^{\xi_{2}},$ then
$\beta\left(  fg\right)  \leq\omega^{\xi},$ where $\xi=\max\left\{  \xi
_{1}+\xi_{2},\,\xi_{2}+\xi_{1}\right\}  .$ These results do not assume the
boundedness of the functions involved.
\end{abstract}
\maketitle

\section{Preliminaries}

Let $K$ be a compact metric space. A \ function $f:K\rightarrow\mathbb{R}$ is
said to be of \emph{Baire class one}, or simply, \emph{Baire-1}, if there
exists a sequence $\left(  f_{n}\right)  $ of real-valued continuous functions
that converges pointwise to $f.$ Let $\frak{B}_{1}\left(  K\right)  $
(respectively, $\mathcal{B}_{1}\left(  K\right)  $) be the set of all
real-valued (respectively, bounded real-valued) Baire-1 functions on $K.$
Several authors have studied Baire-1 functions in terms of ordinal ranks
associated to each function. (See, e.g., \cite{H-O-R}, \cite{K-L} and
\cite{K}). In this paper, we study the relationship between two of these
ordinal ranks, namely the oscillation rank $\beta$ and the convergence rank
$\gamma.$

We begin by recalling the definitions of the indices $\beta$ and $\gamma$.
Suppose that $H$ is a compact metric space, and $f$ is a real-valued function
whose domain contains $H$. For any $\varepsilon>0$, let $H^{0}\left(
f,\varepsilon\right)  =H$. If $H^{\alpha}(f,\varepsilon)$ is defined for some
countable ordinal $\alpha$, let $H^{\alpha+1}\left(  f,\varepsilon\right)  $
be the set of all those $x\in H^{\alpha}\left(  f,\varepsilon\right)  $ such
that for every open set $U$ containing $x,$ there are two points $x_{1}$ and
$x_{2}$ in $U\cap H^{\alpha}\left(  f,\varepsilon\right)  $ with $\left|
f\left(  x_{1}\right)  -f\left(  x_{2}\right)  \right|  \geq\varepsilon.$ For
a countable limit ordinal $\alpha$, we let
\[
H^{\alpha}\left(  f,\varepsilon\right)  =\bigcap_{\alpha^{\prime}<\alpha
}H^{\alpha^{\prime}}\left(  f,\varepsilon\right)  .
\]
The index $\beta_{H}(f,\varepsilon)$ is taken to be $\text{the least }%
\alpha\text{ with }H^{\alpha}\left(  f,\varepsilon\right)  =\emptyset$ if such
$\alpha\text{ exists,}$ and $\omega_{1}$ otherwise. The \textbf{oscillation
index} of $f$ is
\[
\beta_{H}\left(  f\right)  =\sup\left\{  \beta_{H}\left(  f,\varepsilon
\right)  :\varepsilon>0\right\}  .
\]
If the ambient space $H$ is clear from the context, we write $\beta
(f,\varepsilon)$ and $\beta(f)$ in place of $\beta_{H}(f,\varepsilon)$ and
$\beta_{H}(f)$ respectively.

The $\gamma$ index is defined analogously. If $\left(  f_{n}\right)  $ is a
sequence of real-valued functions such that $H\subseteq\bigcap_{n}%
\operatorname*{dom}\left(  f_{n}\right)  ,$ let $H^{0}\left(  \left(
f_{n}\right)  ,\varepsilon\right)  =H$ for any $\varepsilon>0$. If $H^{\alpha
}\left(  \left(  f_{n}\right)  ,\varepsilon\right)  $ has been defined for
some countable ordinal $\alpha$, let $H^{\alpha+1}\left(  \left(
f_{n}\right)  ,\varepsilon\right)  $ be the set of all those $x\in H^{\alpha
}\left(  \left(  f_{n}\right)  ,\varepsilon\right)  $ such that for every open
set $U$ containing $x$ and any $m\in\mathbb{N},$ there are two integers
$n_{1}$,\thinspace\ $n_{2}$ with $n_{1}>n_{2}>m$ and $x^{\prime}\in U\cap
H^{\alpha}\left(  \left(  f_{n}\right)  ,\varepsilon\right)  $ such that
$\left|  f_{n_{1}}\left(  x^{\prime}\right)  -f_{n_{2}}\left(  x^{\prime
}\right)  \right|  \geq\varepsilon.$ Define
\[
H^{\alpha}\left(  \left(  f_{n}\right)  ,\varepsilon\right)  =\bigcap
_{\alpha^{\prime}<\alpha}H^{\alpha^{\prime}}\left(  \left(  f_{n}\right)
,\varepsilon\right)
\]
if $\alpha$ is a countable limit ordinal. Let $\gamma_{H}\left(  \left(
f_{n}\right)  ,\varepsilon\right)  $ be the {least} $\alpha$ with $H^{\alpha
}\left(  \left(  f_{n}\right)  ,\varepsilon\right)  =\emptyset$ {if such}
$\alpha$ {exists,} and $\omega_{1}$ otherwise. Finally, the
\textbf{convergence index} of $\left(  f_{n}\right)  $ is the ordinal
\[
\gamma_{H}\left(  \left(  f_{n}\right)  \right)  =\sup\left\{  \gamma
_{H}\left(  \left(  f_{n}\right)  ,\varepsilon\right)  :\varepsilon>0\right\}
.
\]
Again, if there is no ambiguity about the space $H$, we write $\gamma
((f_{n}),\varepsilon)$ and $\gamma((f_{n}))$ for $\gamma_{H}((f_{n}%
),\varepsilon)$ and $\gamma_{H}((f_{n}))$ respectively.

It is known that a function $f:K\rightarrow\mathbb{R}$ is Baire-1 if and only
if $\beta\left(  f\right)  <\omega_{1}.$ (See \cite[Proposition 1.2]{K-L}.)
Following \cite{K-L}, we define the set of functions of small Baire class
$\xi$ and the set of $\emph{bounded}$ functions of small Baire class $\xi$ for
each countable ordinal $\xi$ as
\[
\frak{B}_{1}^{\xi}\left(  K\right)  =\left\{  f\in\frak{B}_{1}\left(
K\right)  :\beta\left(  f\right)  \leq\omega^{\xi}\right\}
\]
and
\[
\mathcal{B}_{1}^{\xi}\left(  K\right)  =\left\{  f\in\mathcal{B}_{1}\left(
K\right)  :\beta\left(  f\right)  \leq\omega^{\xi}\right\}
\]
respectively. In \cite{K}, the following results are shown.

\begin{theorem}
\label{thmK} Let $K$ be a compact metric space.

\begin{enumerate}
\item \emph{\cite[Theorem 7]{K}} If $\xi$ is a finite ordinal, then a function
$f\in\mathcal{B}_{1}^{\xi+1}\left(  K\right)  $ if and only if there exists a
sequence $(f_{n})$ in $\mathcal{B}_{1}^{1}\left(  K\right)  $ converging
pointwise to $f$ such that $\gamma((f_{n}))\leq\omega^{\xi}$.

\item \emph{\cite[Corollary 9]{K}} If $\xi$ is an infinite countable ordinal,
and $f\in\mathcal{B}_{1}\left(  K\right)  $ is the pointwise limit of a
sequence $(f_{n})$ in $\mathcal{B}_{1}^{1}\left(  K\right)  $ such that
$\gamma((f_{n}))\leq\omega^{\xi}$, then $\beta(f)\leq\omega^{\xi}$.
\end{enumerate}
\end{theorem}

One of our main results generalizes and unifies the two parts of Theorem
\ref{thmK}.

\begin{theorem}
\label{main} Let $K$ be a compact metric space and let $\xi_{1}$, $\xi_{2}$ be
countable ordinals. A function $f\in\frak{B}_{1}^{\xi_{1}+\xi_{2}}\left(
K\right)  ,$ respectively, $\mathcal{B}_{1}^{\xi_{1}+\xi_{2}}\left(  K\right)
,$ if and only if there exists a sequence $(f_{n})$ in $\frak{B}_{1}^{\xi_{1}%
}\left(  K\right)  ,$ respectively, $\mathcal{B}_{1}^{\xi_{1}}\left(
K\right)  ,$ converging pointwise to $f$ such that $\gamma((f_{n}))\leq
\omega^{\xi_{2}}$.
\end{theorem}

In the course of proving Theorem \ref{main}, we show that any Baire-$1$
function $f$ on a closed subspace $H$ of a compact metric space $K$ can be
extended to a Baire-$1$ function $g$ on $K$ such that $\beta_{H}(f)=\beta
_{K}(g)$ (Theorem \ref{EE}). When $\beta_{H}(f)=1$, this is the familiar
Tietze Extension Theorem. Proposition 2.1 and Theorem 2.3 in \cite{K-L} yield
that for a bounded Baire-$1$ function $f$, $\beta(f)$ is the smallest ordinal
$\xi$ such that there exists a sequence of continuous functions $(f_{n})$
converging pointwise to $f$ and having $\gamma((f_{n}))=\xi$. Theorem \ref{18}
below shows that the same result holds without the boundedness assumption on
the function $f$. In the last section, we consider the product of Baire-$1$
functions. In contrast to the class $\mathcal{B}_{1}^{\xi}(K)$, the class
$\frak{B}_{1}^{\xi}\left(  K\right)  $ is not closed under multiplication.
Theorem \ref{Product} shows that if $f\in\frak{B}_{1}^{\xi_{1}}\left(
K\right)  $ and $g\in\frak{B}_{1}^{\xi_{2}}\left(  K\right)  ,$ then
$fg\in\frak{B}_{1}^{\xi}\left(  K\right)  ,$ where $\xi=\max\left\{  \xi
_{1}+\xi_{2},\,\xi_{2}+\xi_{1}\right\}  .$ It is also shown that this result
is the best possible.

Our notation is standard. In the sequel, $K$ will always denote a compact
metric space. If $H$ is a closed subset of $K,$ the derived set $H^{\prime}$
is the set of all limit points of $H.$ A transfinite sequence of derived sets
is defined in the usual manner. Let $H^{\left(  0\right)  }=H$ and $H^{\left(
\alpha+1\right)  }=\left(  H^{\left(  \alpha\right)  }\right)  ^{\prime}$ for
any ordinal $\alpha.$ If $\alpha$ is a limit ordinal, let
\[
H^{\left(  \alpha\right)  }=\bigcap_{\alpha^{\prime}<\alpha}H^{\left(
\alpha^{\prime}\right)  }.
\]
Given real-valued functions $f$ and $g$ defined on a set $S$, we let
\[
\Vert f-g\Vert_{S}=\sup\{|f(s)-g(s)|:s\in S\}.
\]
When there is no cause for confusion, we write $\Vert f-g\Vert$ for $\Vert
f-g\Vert_{S}$. Since we shall be dealing with unbounded functions in general,
this functional can take the value $\infty$ and is not a ``norm''. However, it
is compatible with the topology of uniform convergence on the set
$\mathbb{R}^{S}$ of all real-valued functions on $S$ in the sense that the
sets
\[
U(f,\varepsilon)=\{g:\Vert g-f\Vert_{S}<\varepsilon\}
\]
form a basis for the said topology.

\bigskip

\section{Oscillation and convergence of Baire-1 functions}

We begin by proving a result that yields an upper bound of the oscillation
index of a Baire-1 function $f$ as the product of the convergence index of a
sequence of functions $\left(  f_{n}\right)  $ converging pointwise to $f$,
and the supremum of the oscillation indices of $f_{n}$'s$.$

\begin{lemma}
\label{1}Let $U$ and $L$ be sets such that $U\subseteq L\subseteq K,$ where
$U$ is open in $K$ and $L$ is closed in $K.$ Suppose $f$,$\ f_{n}$ $\left(
n\geq1\right)  $ are Baire-1 functions on $K,$ $\alpha<\omega_{1}$, and
$\varepsilon>0.$ Then

(a) $L^{\alpha}\left(  f,\varepsilon\right)  \subseteq K^{\alpha}\left(
f,\varepsilon\right)  \cap L,$

(b) $L^{\alpha}\left(  \left(  f_{n}\right)  ,\varepsilon\right)  \subseteq
K^{\alpha}\left(  \left(  f_{n}\right)  ,\varepsilon\right)  \cap L,$

(c) $K^{\alpha}\left(  f,\varepsilon\right)  \cap U\subseteq L^{\alpha}\left(
f,\varepsilon\right)  ,$

(d) $K^{\alpha}\left(  \left(  f_{n}\right)  ,\varepsilon\right)  \cap
U\subseteq L^{\alpha}\left(  \left(  f_{n}\right)  ,\varepsilon\right)  .$
\end{lemma}

\begin{proof}
We only prove (c). The proof is by induction on $\alpha.$ The statement is
trivial if $\alpha=0$ or a limit ordinal. Suppose the statement is true for
all ordinals not greater than $\alpha.$ Let $x\in K^{\alpha+1}\left(
f,\varepsilon\right)  \cap U.$ If $N$ is a neighborhood of $x$ in $K$, then
$N\cap U$ is open in $K.$ Thus there exist $x_{1},\,x_{2}\in\left(  N\cap
U\right)  \cap K^{\alpha}\left(  f,\varepsilon\right)  =N\cap\left(  U\cap
K^{\alpha}\left(  f,\varepsilon\right)  \right)  \subseteq N\cap L^{\alpha
}\left(  f,\varepsilon\right)  $ such that $\left|  f\left(  x_{1}\right)
-f\left(  x_{2}\right)  \right|  \geq\varepsilon.$ Hence $x\in L^{\alpha
+1}\left(  f,\varepsilon\right)  .$
\end{proof}

\begin{proposition}
\label{T1}Let $\left(  f_{n}\right)  $ be a sequence in $\frak{B}_{1}\left(
K\right)  $ and let $\varepsilon>0.$ Suppose that $\beta\left(  f_{n}%
,\varepsilon\right)  \leq\beta_{0}$ for all $n\in\mathbb{N}$, and
$\gamma\left(  \left(  f_{n}\right)  ,\varepsilon\right)  \leq\gamma_{0}.$ If
$\left(  f_{n}\right)  $ converges pointwise to a function $f,$ then
$\beta\left(  f,3\varepsilon\right)  \leq\beta_{0}\cdot\gamma_{0}.$
\end{proposition}

\begin{proof}
We first consider the case $\gamma_{0}=1.$ Then $K^{1}\left(  \left(
f_{n}\right)  ,\varepsilon\right)  =\emptyset.$ For each $x\in K,$ there exist
an open neighborhood $U_{x}$ of $x$ and $p_{x}\in\mathbb{N}$ such that
whenever $n>m>p_{x},$%
\[
\left|  f_{n}\left(  x^{\prime}\right)  -f_{m}\left(  x^{\prime}\right)
\right|  <\varepsilon
\]
for all $x^{\prime}\in U_{x}.$ By the compactness of $K,$ there exist
$x_{1},x_{2},...,x_{k}$ such that
\[
K\subseteq\bigcup_{i=1}^{k}U_{x_{i}}.
\]
Let $p_{0}=\max\left\{  p_{x_{1}},p_{x_{2}},...,p_{x_{k}}\right\}  .$ Then for
all $n>m>p_{0}$ and $y\in K,$ we have $y\in U_{x_{i}}$ for some $i,\,1\leq
i\leq k.$ Since $n>m>p_{x_{i}},$
\[
\left|  f_{n}\left(  y\right)  -f_{m}\left(  y\right)  \right|  <\varepsilon.
\]
Taking limit as $n\rightarrow\infty,$ we have
\begin{equation}
\left|  \left|  f-f_{m}\right|  \right|  \leq\varepsilon\text{ \ \ for all
}m>p_{0}. \label{I}%
\end{equation}
Using (\ref{I}), it is easy to verify by induction that
\[
K^{\alpha}\left(  f,3\varepsilon\right)  \subseteq K^{\alpha}\left(
f_{p_{0}+1},\varepsilon\right)
\]
for all $\alpha<\omega_{1}.$ In particular,
\[
K^{\beta_{0}}\left(  f,3\varepsilon\right)  \subseteq K^{\beta_{0}}\left(
f_{p_{0}+1},\varepsilon\right)  =\emptyset.
\]
Hence $\beta\left(  f,3\varepsilon\right)  \leq\beta_{0}=\beta_{0}\cdot
\gamma_{0}.$

Suppose the assertion is true for some $\gamma_{0}.$ Let $\left(
f_{n}\right)  $ be a sequence in $\frak{B}_{1}\left(  K\right)  $ that
converges pointwise to a function $f.$ Suppose there exists $\varepsilon>0$
such that $\beta\left(  f_{n},\varepsilon\right)  \leq\beta_{0}$ for all
$n\in\mathbb{N}$ and $\gamma\left(  \left(  f_{n}\right)  ,\varepsilon\right)
\leq\gamma_{0}+1.$ We need to show $\beta\left(  f,3\varepsilon\right)
\leq\beta_{0}\cdot\left(  \gamma_{0}+1\right)  .$ Since $\gamma\left(  \left(
f_{n}\right)  ,\varepsilon\right)  \leq\gamma_{0}+1,$ we have $K^{\gamma
_{0}+1}\left(  \left(  f_{n}\right)  ,\varepsilon\right)  =\emptyset.$ For
each $m\in\mathbb{N}$, let $U_{m}$ denote the $\frac{1}{m}-$neighborhood of
$K^{\gamma_{0}}\left(  \left(  f_{n}\right)  ,\varepsilon\right)  .$ Denote
$K\setminus U_{m}$ by $\tilde{K}_{m}.$ From Lemma \ref{1}(a) and \ref{1}(b),
for each $n\in\mathbb{N},$ $\beta_{\tilde{K}_{m}}\left(  f_{n},\varepsilon
\right)  \leq\beta_{0}$ and $\gamma_{\tilde{K}_{m}}\left(  \left(
f_{n}\right)  ,\varepsilon\right)  \leq\gamma_{0}$. By the inductive
hypothesis, we see that
\[
\beta_{\tilde{K}_{m}}\left(  f,3\varepsilon\right)  \leq\beta_{0}\cdot
\gamma_{0}.
\]
From this and applying Lemma \ref{1}(c) with $U=K\setminus\overline{U}_{m},$
$L=\tilde{K}_{m}$ for all $m\in\mathbb{N}$, we see that $K^{\beta_{0}%
\cdot\gamma_{0}}\left(  f,3\varepsilon\right)  \subseteq K^{\gamma_{0}}\left(
\left(  f_{n}\right)  ,\varepsilon\right)  .$ Let
\[
\tilde{K}=K^{\beta_{0}\cdot\gamma_{0}}\left(  f,3\varepsilon\right)  \subseteq
K^{\gamma_{0}}\left(  \left(  f_{n}\right)  ,\varepsilon\right)  .
\]
Then
\[
\beta_{\tilde{K}}\left(  f_{n},\varepsilon\right)  \leq\beta_{0}\text{ and
}\gamma_{\tilde{K}}\left(  \left(  f_{n}\right)  ,\varepsilon\right)  =1.
\]
Thus
\[
\beta_{\tilde{K}}\left(  f,3\varepsilon\right)  \leq\beta_{0}\text{ by the
case when }\gamma_{0}=1.
\]
Therefore
\[
K^{\beta_{0}\cdot\left(  \gamma_{0}+1\right)  }\left(  f,3\varepsilon\right)
=K^{\beta_{0}\cdot\gamma_{0}+\beta_{0}}\left(  f,3\varepsilon\right)
=\tilde{K}^{\beta_{0}}\left(  f,3\varepsilon\right)  =\emptyset.
\]
Hence
\[
\beta\left(  f,3\varepsilon\right)  \leq\beta_{0}\cdot\left(  \gamma
_{0}+1\right)  .
\]

Suppose $\gamma_{0}<\omega_{1}$ is a limit ordinal and the statement holds for
all ordinals $\gamma<\gamma_{0}.$ Let $\left(  f_{n}\right)  \subseteq
\frak{B}_{1}\left(  K\right)  $ be a sequence that converges pointwise to a
function $f$ and let $\varepsilon>0$ be given$.$ Suppose that $\beta\left(
f_{n},\varepsilon\right)  \leq\beta_{0}$ for all $n\in\mathbb{N}$, and
$\gamma\left(  \left(  f_{n}\right)  ,\varepsilon\right)  \leq\gamma_{0}.$
Then $\gamma\left(  \left(  f_{n}\right)  ,\varepsilon\right)  <\gamma_{0}$
and $\beta\left(  f,3\varepsilon\right)  \leq\beta_{0}\cdot$ $\gamma\left(
\left(  f_{n}\right)  ,\varepsilon\right)  <\beta_{0}\cdot\gamma_{0}.$
\end{proof}

\begin{theorem}
\label{4}Let $\left(  f_{n}\right)  $ be a sequence $\frak{B}_{1}\left(
K\right)  $ converging pointwise to a function $f.$ Suppose $\sup\left\{
\beta\left(  f_{n}\right)  :n\in\mathbb{N}\right\}  \leq\beta_{0}$ and
$\gamma\left(  \left(  f_{n}\right)  \right)  \leq\gamma_{0}.$ Then $f$ is
Baire-1 and $\beta\left(  f\right)  \leq\beta_{0}\cdot\gamma_{0}.$
\end{theorem}

For the next corollary, recall that $DBSC\left(  K\right)  $ is the space of
all differences of semicontinuous functions on $K.$ It is known that
$\mathcal{B}_{1}^{1}\left(  K\right)  $ is the closure of $DBSC\left(
K\right)  $ in the topology of uniform convergence (\cite[Theorem 3.1]{K-L}).

\begin{corollary}
[{\cite[Corollary 9]{K}}]Let $f\in\mathcal{B}_{1}\left(  K\right)  $ be the
pointwise limit of a sequence $\left(  f_{n}\right)  \subseteq DBSC\left(
K\right)  $. If $\gamma\left(  \left(  f_{n}\right)  \right)  \leq\omega^{\xi
},$ $\omega\leq\xi<\omega_{1},$ then $\beta\left(  f\right)  \leq\omega^{\xi}.$
\end{corollary}

\section{Extension of Baire-1 functions}

In this section, we establish several results regarding the extension of
Baire-1 functions. They are analogs of the Tietze Extension Theorem for
continuous functions. These results are applied in the next section in proving
the converse of Theorem \ref{4}.

\begin{lemma}
\label{5}Suppose that $F$ is a closed subspace of $K$ and that $f$ is a
Baire-1 function on $F$. For any $\varepsilon>0$, there exists a continuous
function $g:K\setminus F^{1}\left(  f,\varepsilon\right)  \rightarrow
\mathbb{R}$ such that
\[
\left\|  g-f\right\|  _{F\setminus F^{1}\left(  f,\varepsilon\right)  }%
\leq\varepsilon.
\]
\end{lemma}

\begin{proof}
For any $x\in F\setminus F^{1}\left(  f,\varepsilon\right)  ,$ choose an open
neighborhood $U_{x}$ of $x$ in $K$ such that $U_{x}\cap F^{1}\left(
f,\varepsilon\right)  =\emptyset$ and $\left|  f\left(  x_{1}\right)
-f\left(  x_{2}\right)  \right|  <\varepsilon\text{ for all }x_{1},\,x_{2}\in
U_{x}\cap F.$ The collection $\mathcal{U}=\left\{  U_{x}:x\in F\setminus
F^{1}\left(  f,\varepsilon\right)  \right\}  \cup\left\{  K\setminus
F\right\}  $ is an open cover of $K\setminus F^{1}\left(  f,\varepsilon
\right)  $. By \cite{D}, Theorems IX.5.3 and VIII.4.2, there exists a
partition of unity $\left(  \varphi_{U}\right)  _{U\in\mathcal{U}}$
subordinated to $\mathcal{U}.$ If $U=U_{x}\in\mathcal{U}$ for some $x\in
F\setminus F^{1}\left(  f,\varepsilon\right)  $, let $a_{U}=f(x);$ if
$U=K\setminus F$, let $a_{U}=0$. Define $g:K\setminus F^{1}\left(
f,\varepsilon\right)  \rightarrow\mathbb{R}$ by $g=\sum_{U\in\mathcal{U}}%
a_{U}\varphi_{U}.$ The sum is well-defined since $\left\{
\text{supp\thinspace}\varphi_{U}:U\in\mathcal{U}\right\}  $ is locally finite.
Let $x\in F\setminus F^{1}\left(  f,\varepsilon\right)  .$ Then $\mathcal{V}%
=\left\{  U\in\mathcal{U}:\varphi_{U}\left(  x\right)  \neq0\right\}  $ is a
finite set, $\varphi_{U}\left(  x\right)  >0$ for all $U\in\mathcal{V}$ and
$\sum_{U\in\mathcal{V}}\varphi_{U}\left(  x\right)  =1.$ If $U\in\mathcal{V},$
then $x\in U\cap F;$ hence $U\neq K\setminus F$. Therefore, $U=U_{y}$ for some
$y\in F\setminus F^{1}\left(  f,\varepsilon\right)  .$ But then $x,y\in
U_{y}\cap F$ implies that $\left|  a_{U}-f\left(  x\right)  \right|  =\left|
f\left(  y\right)  -f\left(  x\right)  \right|  <\varepsilon.$ It follows
that
\begin{align*}
\left|  g\left(  x\right)  -f\left(  x\right)  \right|   &  =\left|
\sum_{U\in\mathcal{U}}a_{U}\varphi_{U}\left(  x\right)  -f\left(  x\right)
\right|  =\left|  \sum_{U\in\mathcal{V}}a_{U}\varphi_{U}\left(  x\right)
-\sum_{U\in\mathcal{V}}f\left(  x\right)  \varphi_{U}\left(  x\right)  \right|
\\
&  \leq\sum_{U\in\mathcal{V}}\left|  a_{U}-f\left(  x\right)  \right|
\varphi_{U}\left(  x\right)  <\varepsilon.
\end{align*}
This shows that
\[
\left\|  g-f\right\|  _{F\setminus F^{1}\left(  f,\varepsilon\right)  }%
\leq\varepsilon.
\]
Finally, if $x$ is a point in $K\setminus F^{1}\left(  f,\varepsilon\right)
$, there exists an open neighborhood $V$ of $x$ in $K$ such that $V\cap
F^{1}\left(  f,\varepsilon\right)  =\emptyset$ and $\mathcal{W}=\left\{
U\in\mathcal{U}:\text{supp\thinspace}\varphi_{U}\cap V\neq\emptyset\right\}  $
is finite. Now
\[
g_{|V}=\sum_{U\in\mathcal{U}}a_{U}\varphi_{U\,|V}=\sum_{U\in\mathcal{W}}%
a_{U}\varphi_{U\,|V}.
\]
Hence $g_{|V}$ is continuous on $V$, since it is a finite linear combination
of continuous functions. In particular, $g$ is continuous at $x$. As $x\in
K\setminus F^{1}\left(  f,\varepsilon\right)  $ is arbitrary, $g$ is
continuous on $K\setminus F^{1}\left(  f,\varepsilon\right)  $.
\end{proof}

\begin{theorem}
\label{AE}Suppose that $F$ is a closed subspace of $K$ and that $f$ is a
Baire-1 function on $F$. For any $1\leq\beta_{0}<\omega_{1},$ and any
$\varepsilon>0,$ there exists $g:K\setminus F^{\beta_{0}}\left(
f,\varepsilon\right)  \rightarrow\mathbb{R}$ such that
\[
\left\|  g-f\right\|  _{F\setminus F^{\beta_{0}}\left(  f,\varepsilon\right)
}\leq\varepsilon
\]
and
\[
\beta_{H}\left(  g\right)  \leq\beta_{0}\text{ for all compact subsets
}H\text{ of }K\setminus F^{\beta_{0}}\left(  f,\varepsilon\right)  .
\]
\end{theorem}

\begin{proof}
Let $h:K\setminus F^{1}\left(  f,\varepsilon\right)  \rightarrow\mathbb{R}$ be
the function obtained from Lemma \ref{5}. If $1\leq\alpha<\beta_{0}$, let
$\tilde{K}=\tilde{F}=F^{\alpha}\left(  f,\varepsilon\right)  $. Applying Lemma
\ref{5} with $\tilde{K}$, $\tilde{F}$, and the function $f$ yields a
continuous function $g_{\alpha}:F^{\alpha}\left(  f,\varepsilon\right)
\setminus F^{\alpha+1}\left(  f,\varepsilon\right)  \rightarrow\mathbb{R}$
such that
\[
\left\|  g_{\alpha}-f\right\|  _{F^{\alpha}\left(  f,\varepsilon\right)
\setminus F^{\alpha+1}\left(  f,\varepsilon\right)  }\leq\varepsilon.
\]
Let $g=h\cup\left(  \bigcup_{\alpha<\beta_{0}}g_{\alpha}\right)  :K\setminus
F^{\beta_{0}}\left(  f,\varepsilon\right)  \rightarrow\mathbb{R}$. Then
$\left\|  g-f\right\|  _{F\setminus F^{\beta_{0}}\left(  f,\varepsilon\right)
}\leq\varepsilon.$

Suppose that $\delta>0$ and $H$ is a compact subset of $K\setminus
F^{\beta_{0}}\left(  f,\varepsilon\right)  .$ If $x\notin F^{1}\left(
f,\varepsilon\right)  ,$ then there exists an open neighborhood $U$ of $x$
such that
\[
\overline{U}\cap F^{1}\left(  f,\varepsilon\right)  =\emptyset.
\]
Note that $g_{|\overline{U}}=h_{|\overline{U}}.$ By Lemma \ref{1}(c),%
\[
H^{1}\left(  g,\delta\right)  \cap U\subseteq\left(  H\cap\overline{U}\right)
^{1}\left(  g,\delta\right)  =\left(  H\cap\overline{U}\right)  ^{1}\left(
h,\delta\right)  =\emptyset
\]
by the continuity of $h.$ In particular, $x\notin H^{1}\left(  g,\delta
\right)  .$ It follows that
\[
H^{1}\left(  g,\delta\right)  \subseteq H\cap F^{1}\left(  f,\varepsilon
\right)  .
\]
Repeating the argument inductively yields that%
\[
H^{\beta_{0}}\left(  g,\delta\right)  \subseteq H\cap F^{\beta_{0}}\left(
f,\varepsilon\right)  =\emptyset.
\]
Hence $\beta_{H}\left(  g\right)  \leq\beta_{0}$, as required.
\end{proof}

\bigskip

We obtain the following corollaries by taking $F=K$ and $\beta_{0}=\beta
_{F}\left(  f\right)  $ respectively.

\begin{corollary}
\label{c9}Let $f$ be a Baire-1 function on $K$ such that $\beta\left(
f,\varepsilon\right)  \leq\beta_{0}$ for some $1\leq\beta_{0}<\omega_{1}$ and
$\varepsilon>0.$ Then there exists $g:K\rightarrow\mathbb{R}$ such that
\[
\left\|  g-f\right\|  \leq\varepsilon\text{ and }\beta\left(  g\right)
\leq\beta_{0}.
\]
\end{corollary}

\begin{corollary}
\label{C4}Let $F$ be a closed subspace of $K.$ If $f$ is a Baire-1 function on
$F$, then for every $\varepsilon>0$ there exists a Baire-1 function $g$ on $K$
such that
\[
\left\|  g-f\right\|  _{F}\leq\varepsilon\text{ and }\beta_{K}\left(
g\right)  \leq\beta_{F}\left(  f\right)  .
\]
\end{corollary}

Next we show that Corollary \ref{C4} can be improved to an exact extension
theorem (i.e., the case $\varepsilon=0$). In the statement of Lemma \ref{12},
the vacuous sum $\sum_{j=1}^{0}g_{j}$ is taken to be the zero function.

\begin{lemma}
\label{12}Let $F$ be a closed subspace of $K$ and let $f$ be a Baire-1
function on $F.$ Then there exists a sequence of Baire-1 functions $\left(
g_{n}\right)  $ on $K$ such that

(a) $g_{n}$ is continuous on $K\setminus F^{1}\left(  f-\sum_{j=1}^{n-1}%
g_{j},\frac{1}{2^{n-1}}\right)  $ for all $n\in\mathbb{N},$

(b) $\left\|  f-\sum_{j=1}^{n}g_{j}\right\|  _{F\setminus F^{1}\left(
f,\frac{1}{4^{n-1}}\right)  }\leq\dfrac{1}{2^{n-1}},$ $n\in\mathbb{N},$

(c) $\left\|  g_{n}\right\|  _{K}\leq\dfrac{1}{2^{n-2}}$ if $n\geq2,$ and

(d) $F^{1}\left(  f-\sum_{j=1}^{n}g_{j},\delta\right)  \subseteq F^{1}\left(
f,\frac{\delta}{2^{n}}\right)  $ if $0<\delta\leq\frac{1}{2^{n-2}},$
$n\in\mathbb{N}.$
\end{lemma}

\begin{proof}
The functions $\left(  g_{n}\right)  $ are constructed inductively$.$ By Lemma
\ref{5}, there exists a continuous function $g_{1}:K\setminus F^{1}\left(
f,1\right)  \rightarrow\mathbb{R}$ such that $\left\|  f-g_{1}\right\|
_{F\setminus F^{1}\left(  f,1\right)  }\leq1$. Extend $g_{1}$ to a function on
$K$ by defining $g_{1}$ to be $0$ on $F^{1}\left(  f,1\right)  $. Then (a) and
(b) hold. Condition (c) holds vacuously. Moreover, if $x\in F\setminus
F^{1}\left(  f,\frac{\delta}{2}\right)  ,\,0<\delta\leq2,$ then there exists a
neighborhood $U_{1}$ of $x$ in $F$ such that $\left|  f\left(  x_{1}\right)
-f\left(  x_{2}\right)  \right|  <\frac{\delta}{2}$ for all $x_{1},$ $x_{2}\in
U_{1}.$ Note that since $x\in F\setminus F^{1}\left(  f,\frac{\delta}%
{2}\right)  ,$ $g_{1}$ is continuous at $x.$ Hence there exists a neighborhood
$U_{2}$ of $x$ in $F$ such that $\left|  g_{1}\left(  x_{1}\right)
-g_{1}\left(  x_{2}\right)  \right|  <\frac{\delta}{2}$ for all $x_{1},$
$x_{2}\in U_{2}.$ Let $U=U_{1}\cap U_{2}.$ Then $U$ is a neighborhood of $x$
in $F.$ For all $x_{1},\,x_{2}\in U,$%
\[
\left|  \left(  f-g_{1}\right)  \left(  x_{1}\right)  -\left(  f-g_{1}\right)
\left(  x_{2}\right)  \right|  <\delta.
\]
Hence $x\notin F\left(  f-g_{1},\delta\right)  .$ This proves (d).

Suppose that $g_{1},\,g_{2},\,...,\,g_{n}$ have been chosen. By Lemma \ref{5},
there exists a continuous function $h:K\setminus F^{1}\left(  f-\sum_{j=1}%
^{n}g_{j},\frac{1}{2^{n}}\right)  \rightarrow\mathbb{R}$ such that
\[
\left\|  f-\sum_{j=1}^{n}g_{j}-h\right\|  _{F\setminus F^{1}(f-%
{\textstyle\sum_{j=1}^{n}}
g_{j},\frac{1}{2^{n}})}\leq\frac{1}{2^{n}}.
\]
Define $\tilde{h}$ on $K\setminus F^{1}(f-%
{\textstyle\sum_{j=1}^{n}}
g_{j},\frac{1}{2^{n}})$ by $\tilde{h}=\left(  h\wedge\dfrac{1}{2^{n-1}%
}\right)  \vee\dfrac{-1}{2^{n-1}}.$ Then $\tilde{h}$ is continuous on
$K\setminus F^{1}(f-%
{\textstyle\sum_{j=1}^{n}}
g_{j},\frac{1}{2^{n}}).$ By (d), $F^{1}(f-%
{\textstyle\sum_{j=1}^{n}}
g_{j},\frac{1}{2^{n}})\subseteq F^{1}\left(  f,\dfrac{1}{4^{n}}\right)  .$
Hence $\tilde{h}$ is defined and continuous on $K\setminus F^{1}\left(
f,\dfrac{1}{4^{n}}\right)  .$ Moreover, it follows from (b) that
\begin{equation}
\left\|  f-\sum_{j=1}^{n}g_{j}\right\|  _{F\setminus F^{1}\left(  f,\frac
{1}{4^{n}}\right)  }\leq\dfrac{1}{2^{n-1}}. \label{2}%
\end{equation}
From inequality (\ref{2}) and the definition of $\tilde{h},$ we have
\[
\left\|  f-\sum_{j=1}^{n}g_{j}-\tilde{h}\right\|  _{F\setminus F^{1}\left(
f,\frac{1}{4^{n}}\right)  }\leq\left\|  f-\sum_{j=1}^{n}g_{j}-h\right\|
_{F\setminus F^{1}\left(  f,\frac{1}{4^{n}}\right)  }.
\]
Therefore, $\left\|  f-\sum_{j=1}^{n}g_{j}-\tilde{h}\right\|  _{F\setminus
F^{1}\left(  f,\frac{1}{4^{n}}\right)  }\leq\dfrac{1}{2^{n}}.$ Now define
\[
g_{n+1}=\left\{
\begin{array}
[c]{cc}%
\tilde{h} & \text{on }K\setminus F^{1}(f-%
{\textstyle\sum_{j=1}^{n}}
g_{j},\frac{1}{2^{n}})\\
0 & \text{otherwise}%
\end{array}
.\right.
\]
Then $g_{n+1}$ is continuous on $K\setminus F^{1}(f-%
{\textstyle\sum_{j=1}^{n}}
g_{j},\frac{1}{2^{n}}).$ This proves (a). Furthermore,
\[
\left\|  f-\sum_{j=1}^{n+1}g_{j}\right\|  _{F\setminus F^{1}\left(  f,\frac
{1}{4^{n}}\right)  }=\left\|  f-\sum_{j=1}^{n}g_{j}-\tilde{h}\right\|
_{F\setminus F^{1}\left(  f,\frac{1}{4^{n}}\right)  }\leq\dfrac{1}{2^{n}}.
\]
This proves (b). Also,
\[
\left\|  g_{n+1}\right\|  _{K}\leq\left\|  \tilde{h}\right\|  _{K\setminus
F^{1}(f-%
{\textstyle\sum_{j=1}^{n}}
g_{j},\frac{1}{2^{n}})}\leq\dfrac{1}{2^{n-1}}%
\]
by the definition of $\tilde{h}.$ This proves (c). Finally, suppose
$0<\delta\leq\dfrac{1}{2^{n-1}}.$ Assume that $x\in F\setminus F^{1}\left(
f,\dfrac{\delta}{2^{n+1}}\right)  .$ Then $x\notin F^{1}\left(  f-\sum
_{j=1}^{n}g_{j},\dfrac{\delta}{2}\right)  .$ Thus there exists a neighborhood
$U_{1}$ of $x$ in\thinspace$F$ such that
\[
\left|  \left(  f-\sum_{j=1}^{n}g_{j}\right)  \left(  x_{1}\right)  -\left(
f-\sum_{j=1}^{n}g_{j}\right)  \left(  x_{2}\right)  \right|  <\dfrac{\delta
}{2}%
\]
whenever $x_{1},\,x_{2}\in U_{1}.$ Note that since $x\in F\setminus
F^{1}\left(  f-\sum_{j=1}^{n}g_{j},\dfrac{\delta}{2}\right)  ,$ $g_{n+1}$ is
continuous at $x.$ Therefore, there exists a neighborhood $U_{2}$ of $x$ in
$F$ such that $\left|  g_{n+1}\left(  x_{1}\right)  -g_{n+1}\left(
x_{2}\right)  \right|  <\dfrac{\delta}{2}$ for all $x_{1},\,x_{2}\in U_{2}.$
Let $U=U_{1}\cap U_{2}.$ Then $U$ is a neighborhood of $x$ in $F$ such that
\[
\left|  \left(  f-\sum_{j=1}^{n+1}g_{j}\right)  \left(  x_{1}\right)  -\left(
f-\sum_{j=1}^{n+1}g_{j}\right)  \left(  x_{2}\right)  \right|  <\delta
\]
whenever $x_{1},\,x_{2}\in U.$ Hence $x\notin F^{1}\left(  f-\sum_{j=1}%
^{n+1}g_{j},\delta\right)  .$ This proves (d).
\end{proof}

\begin{theorem}
\label{EE}Let $F$ be a closed subspace of $K$ and let $f$ be a Baire-1
function on $F.$ Then there exists a Baire-1 function $g$ on $K$ such that
\[
g_{|F}=f\text{ and }\beta\left(  g\right)  =\beta_{F}\left(  f\right)  .
\]
\end{theorem}

\begin{proof}
Let $\left(  g_{n}\right)  $ be the sequence given by Lemma \ref{12}. Define
$g$ on $K$ by
\[
g=\left\{
\begin{array}
[c]{cc}%
\sum_{j=1}^{\infty}g_{j} & \text{on }K\setminus F\\
f & \text{on }F
\end{array}
\right.  .
\]
Note that by (c) of Lemma \ref{12}, $\sum_{j=1}^{\infty}g_{j}$ converges
uniformly on $K.$ Hence $g$ is well defined. Obviously, $g_{|F}=f.$

\noindent\textbf{Claim.} $K^{1}\left(  g,\frac{1}{2^{n-3}}\right)  \subseteq
F^{1}\left(  f,\frac{1}{4^{n}}\right)  $ for all $n\in\mathbb{N}.$

\noindent\textit{Proof of Claim.}\textsl{ }Let $x\in K\setminus F^{1}\left(
f,\frac{1}{4^{n}}\right)  .$ We consider two cases. Suppose $x\notin F.$ By
Lemma \ref{12}(a), $g_{j}$ is continuous on $K\setminus F$ for all $j.$ Since
$\sum_{j=1}^{n}g_{j}$ converges uniformly to $g$ on $K\setminus F,$ and
$K\setminus F$ is an open subset of $K,$ $g$ is continuous at $x.$ Hence
$x\notin K^{1}\left(  g,\frac{1}{2^{n-3}}\right)  .$ Now suppose $x\in F.$
Then $x\in F\setminus F^{1}\left(  f,\frac{1}{4^{n}}\right)  .$ There is a
neighborhood $U_{1}$ of $x$ in $K$ such that $\left|  f\left(  x\right)
-f\left(  x^{\prime}\right)  \right|  <\dfrac{1}{4^{n}}$ for all $x^{\prime
}\in U_{1}\cap F.$ Also, for $1\leq k\leq n,$
\begin{align*}
F^{1}\left(  f-\sum_{j=1}^{k}g_{j},\dfrac{1}{2^{k}}\right)   &  \subseteq
F^{1}\left(  f,\frac{1}{4^{k}}\right)  \text{ by Lemma \ref{12}(d),}\\
&  \subseteq F^{1}\left(  f,\frac{1}{4^{n}}\right)  .
\end{align*}
Since $g_{k+1}$ is continuous on $K\setminus F^{1}(f-%
{\textstyle\sum_{j=1}^{k}}
g_{j},\frac{1}{2^{k}}),$ $g_{k+1}$ is continuous on $K\setminus F^{1}\left(
f,\frac{1}{4^{n}}\right)  $ for all $k,$ $1\leq k\leq n.$ Similarly,
$F^{1}\left(  f,1\right)  \subseteq F^{1}\left(  f,\frac{1}{4^{n}}\right)  $
and $g_{1}$ is continuous on $K\setminus F^{1}\left(  f,1\right)  $ by Lemma
\ref{12}(a); thus, $g_{1}$ is continuous on $K\setminus F^{1}\left(
f,\frac{1}{4^{n}}\right)  .$ Hence there exists a neighborhood $U_{2}$ of $x$
in $K$ such that $U_{2}\subseteq K\setminus F^{1}\left(  f,\frac{1}{4^{n}%
}\right)  $ and
\[
\left|  \sum_{j=1}^{n+1}g_{j}\left(  x^{\prime}\right)  -\sum_{j=1}^{n+1}%
g_{j}\left(  x\right)  \right|  <\dfrac{1}{2^{n}}\text{ for all }x^{\prime}\in
U_{2}.
\]
Let $U=U_{1}\cap U_{2}.$ Then $U$ is a neighborhood of $x$ in $K$. If
$x^{\prime}\in U\cap F,$ then $x^{\prime}\in U_{1}\cap F.$ Thus $\left|
g\left(  x^{\prime}\right)  -g\left(  x\right)  \right|  =\left|  f\left(
x^{\prime}\right)  -f\left(  x\right)  \right|  <\dfrac{1}{4^{n}}<\dfrac
{1}{2^{n-2}}.$ If $x^{\prime}\in U\setminus F,$ then
\begin{align*}
\left|  g\left(  x^{\prime}\right)  -g\left(  x\right)  \right|   &  =\left|
\sum_{j=1}^{\infty}g_{j}\left(  x^{\prime}\right)  -f\left(  x\right)  \right|
\\
&  \leq\left|  \sum_{j=1}^{n+1}g_{j}\left(  x^{\prime}\right)  -\sum
_{j=1}^{n+1}g_{j}\left(  x\right)  \right|  +\left|  \sum_{j=1}^{n+1}%
g_{j}\left(  x\right)  -f\left(  x\right)  \right|  +\left|  \sum
_{j=n+2}^{\infty}g_{j}\left(  x^{\prime}\right)  \right| \\
&  <\dfrac{1}{2^{n}}+\left|  \sum_{j=1}^{n+1}g_{j}\left(  x\right)  -f\left(
x\right)  \right|  +\sum_{j=n+2}^{\infty}\left\|  g_{j}\right\|  \text{ since
}x^{\prime}\in U_{2},\\
&  \leq\dfrac{1}{2^{n}}+\left|  \sum_{j=1}^{n+1}g_{j}\left(  x\right)
-f\left(  x\right)  \right|  +\sum_{j=n+2}^{\infty}\dfrac{1}{2^{j-2}},\text{
by Lemma \ref{12}(c),}\\
&  \leq\dfrac{1}{2^{n}}+\dfrac{1}{2^{n}}+\dfrac{1}{2^{n-1}},\text{ by Lemma
\ref{12}(b), since }x\in F\setminus F^{1}\left(  f,\frac{1}{4^{n}}\right)  ,\\
&  =\dfrac{1}{2^{n-2}}.
\end{align*}
Thus $\left|  g\left(  x^{\prime}\right)  -g\left(  x\right)  \right|
<\dfrac{1}{2^{n-2}}$ if $x^{\prime}\in U.$ Hence $\left|  g\left(
x_{1}\right)  -g\left(  x_{2}\right)  \right|  <\dfrac{1}{2^{n-3}}$ whenever
$x_{1},$ $x_{2}\in U.$ Therefore $x\notin K^{1}\left(  g,\dfrac{1}{2^{n-3}%
}\right)  .$ This proves the claim.

It follows by induction that $K^{\alpha}\left(  g,\dfrac{1}{2^{n-3}}\right)
\subseteq F^{\alpha}\left(  f,\dfrac{1}{4^{n}}\right)  $ for $1\leq
\alpha<\omega_{1}.$ Indeed, the Claim yields the assertion for $\alpha=1.$ If
the inclusion holds for some $\alpha,$ $1\leq\alpha<\omega_{1},$ let
$\tilde{F}=F^{\alpha}\left(  f,\dfrac{1}{4^{n}}\right)  .$ Then $K^{\alpha
+1}\left(  g,\dfrac{1}{2^{n-3}}\right)  \subseteq\tilde{F}^{1}\left(
g,\dfrac{1}{2^{n-3}}\right)  =\tilde{F}^{1}\left(  f,\dfrac{1}{2^{n-3}%
}\right)  \subseteq\tilde{F}^{1}\left(  f,\dfrac{1}{4^{n}}\right)
=F^{\alpha+1}\left(  f,\dfrac{1}{4^{n}}\right)  .$ Hence the inclusion holds
for $\alpha+1.$ If the inclusion holds for all $1\leq\alpha^{\prime}<\alpha,$
where $\alpha<\omega_{1}$ is a limit ordinal, then
\[
K^{\alpha}\left(  g,\dfrac{1}{2^{n-3}}\right)  =\bigcap_{1\leq\alpha^{\prime
}<\alpha}K^{\alpha^{\prime}}\left(  g,\dfrac{1}{2^{n-3}}\right)
\subseteq\bigcap_{1\leq\alpha^{\prime}<\alpha}F^{\alpha^{\prime}}\left(
f,\dfrac{1}{4^{n}}\right)  =F^{\alpha}\left(  f,\dfrac{1}{4^{n}}\right)  .
\]
This proves the inclusion for $1\leq\alpha<\omega_{1}.$ In particular, if
$\beta_{F}\left(  f\right)  =\beta_{0},$ then $K^{\beta_{0}}\left(
g,\dfrac{1}{2^{n-3}}\right)  \subseteq F^{\beta_{0}}\left(  f,\dfrac{1}{4^{n}%
}\right)  =\emptyset.$ Thus $\beta_{K}\left(  g,\dfrac{1}{2^{n-3}}\right)
\leq\beta_{0}$ for all $n\in\mathbb{N}.$ Hence $\beta_{K}\left(  g\right)
\leq\beta_{0}.$ Of course, since $g_{|F}=f,$ $\beta_{K}\left(  g\right)
\geq\beta_{F}\left(  f\right)  \geq\beta_{0}.$ Therefore $\beta_{K}\left(
g\right)  =\beta_{0}=\beta_{F}\left(  f\right)  .$
\end{proof}

\begin{remark}
\emph{If }$\beta_{F}\left(  f\right)  =1,$\emph{ Theorem \ref{EE} is the
familiar Tietze Extension Theorem. If }$\beta_{F}\left(  f\right)  $\emph{ is
transfinite, the conclusion of Theorem \ref{EE} can be obtained easily by
defining the extension }$g$\emph{ to be }$0$\emph{ on }$K\setminus F.$\emph{
However, we do not see a simple proof for finite }$\beta_{F}\left(  f\right)  .$
\end{remark}

\section{Decomposition of Baire-1 functions}

In this section, we give a proof of Theorem \ref{main}. The extension results
in \S3 are employed in the course of the proof.

\begin{theorem}
Let $f$ be a Baire-1 function on $K,$ $1\leq\beta_{0},\,\gamma_{0}<\omega_{1}$
and $\varepsilon>0.$ Then there exist
\[
\tilde{f}:K\setminus K^{\beta_{0}\cdot\gamma_{0}}\left(  f,\varepsilon\right)
\rightarrow\mathbb{R}%
\]
and
\[
f_{n}:K\setminus K^{\beta_{0}\cdot\gamma_{0}}\left(  f,\varepsilon\right)
\rightarrow\mathbb{R}%
\]
such that $(f_{n})$ converges to ${f}$ pointwise$,$ $\left\|  \tilde
{f}-f\right\|  _{K\setminus K^{\beta_{0}\cdot\gamma_{0}}\left(  f,\varepsilon
\right)  }\leq\varepsilon$ and $\beta_{H}\left(  f_{n}\right)  \leq\beta_{0},$
$\gamma_{H}\left(  \left(  f_{n}\right)  \right)  \leq\gamma_{0}$ for all
compact subsets $H$ of $K\setminus K^{\beta_{0}\cdot\gamma_{0}}\left(
f,\varepsilon\right)  .$
\end{theorem}

\begin{proof}
For $\alpha\leq\gamma_{0},$ let $K_{\alpha}=K^{\beta_{0}\cdot\alpha}\left(
f,\varepsilon\right)  .$ If $n\in\mathbb{N},$ let $U_{n}^{\alpha}$ be the
$\dfrac{1}{n}-$neighborhood of $K_{\alpha}$ in $K.$ For $\alpha<\gamma_{0},$
it follows from Theorem \ref{AE} that there exists $g_{\alpha}:K_{\alpha
}\setminus K_{\alpha+1}\rightarrow\mathbb{R}$ such that $\left\|  g_{\alpha
}-f\right\|  _{K_{\alpha}\setminus K_{\alpha+1}}\leq\varepsilon$ and
$\beta_{H}\left(  g_{\alpha}\right)  \leq\beta_{0}$ for all compact subsets
$H$ of $K_{\alpha}\setminus K_{\alpha+1}.$ List the ordinals in $[0,\gamma
_{0})$ in a (possibly finite) sequence $\left(  \alpha_{n}\right)  _{n=1}%
^{p}.$ Here $p\in\mathbb{N}$ or $p=\infty.$ For each $n\in\mathbb{N},$ let
$F_{n}=\bigcup_{j=1}^{n\wedge p}\left(  K_{\alpha_{j}}\setminus U_{n}%
^{\alpha_{j}+1}\right)  .$ Then $F_{n}$ is a closed subset of $K.$ It is also
easy to see that $K_{\alpha}\setminus U_{n}^{\alpha+1}$ and $K_{\alpha
^{\prime}}\setminus U_{n}^{\alpha^{\prime}+1}$ are disjoint if $\alpha
\neq\alpha^{\prime}.$ Thus $\left(  K_{\alpha_{j}}\setminus U_{n}^{\alpha
_{j}+1}\right)  _{j=1}^{n\wedge p}$ is a partition of $F_{n}$ into clopen (in
$F_{n}$) subsets. Now define $\tilde{g}_{n}:F_{n}\rightarrow K$ to be
$\bigcup_{j=1}^{n\wedge p}g_{\alpha_{j}|K_{\alpha_{j}}\setminus U_{n}%
^{\alpha_{j}+1}}.$ Since $H=K_{\alpha_{j}}\setminus U_{n}^{\alpha_{j}+1}$ is a
compact subset of $K_{\alpha_{j}}\setminus K_{\alpha_{j}+1},$ $\beta
_{H}\left(  g_{\alpha_{j}}\right)  \leq\beta_{0}.$ From the clopeness of the
partition $\left(  K_{\alpha_{j}}\setminus U_{n}^{\alpha_{j}+1}\right)
_{j=1}^{n\wedge p},$ it follows readily that $\beta_{F_{n}}\left(  \tilde
{g}_{n}\right)  \leq\beta_{0}.$ By Theorem \ref{EE}, there exists a function
$f_{n}^{\prime}$ on $K$ such that $f_{n|F_{n}}^{\prime}=\tilde{g}_{n}$ and
$\beta_{K}\left(  f_{n}^{\prime}\right)  \leq\beta_{0}.$ Finally, define
$f_{n}$ to be $f_{n|K\setminus K_{\gamma_{0}}}^{\prime}$ and $\tilde{f}$ to be
$\bigcup_{\alpha<\gamma_{0}}g_{\alpha|K_{\alpha}\setminus K_{\alpha+1}}.$ It
follows from the choices of the $g_{\alpha}$'s that $\left\|  f-\tilde
{f}\right\|  _{K\setminus K_{\gamma_{0}}}\leq\varepsilon.$ Since
$\bigcup_{n=1}^{\infty}F_{n}=K\setminus K_{\gamma_{0}},$ $\lim f_{n}=\tilde
{f}$ pointwise on $K\setminus K_{\gamma_{0}}.$ Suppose $H$ is a compact subset
of $K\setminus K_{\gamma_{0}}.$ Then $\beta_{H}\left(  f_{n}\right)  \leq
\beta_{K}\left(  f_{n}^{\prime}\right)  \leq\beta_{0}.$ To complete the proof,
we claim that for any $\delta>0$ and any $\gamma\leq\gamma_{0},$ $H^{\gamma
}\left(  \left(  f_{n}\right)  ,\delta\right)  \subseteq K_{\gamma}.$ The
proof of this is by induction on $\gamma.$ The case $\gamma=0$ and the limit
case is trivial. Now assume that the claim holds for some $\gamma<\gamma_{0}.$
Let $x\in H^{\gamma}\left(  \left(  f_{n}\right)  ,\delta\right)  \setminus
K_{\gamma+1}.$ Choose $j_{1},\,j_{2}\in\mathbb{N}$ such that $\alpha_{j_{1}%
}=\gamma$ and $d\left(  x,K_{\gamma+1}\right)  \geq\dfrac{1}{j_{2}},$ where
$d$ is the metric on $K.$ Denote $H^{\gamma}\left(  \left(  f_{n}\right)
,\delta\right)  $ by $L$ and the $\dfrac{1}{2j_{0}}$-ball in $K$ centered at
$x$ by $U,$ where $j_{0}=\max\left\{  j_{1},2j_{2}\right\}  .$ Note that
$L\subseteq K_{\gamma}$ by the inductive hypothesis: For all $n\geq j_{0}%
=\max\left\{  j_{1},2j_{2}\right\}  ,$%
\[
L\cap U\subseteq L\cap\overline{U}\subseteq K_{\alpha_{j_{1}}}\setminus
U_{n}^{\alpha_{j_{1}}+1}\subseteq F_{n}.
\]
This implies that $f_{n|L\cap\overline{U}}=\tilde{g}_{n|L\cap\overline{U}%
}=g_{\alpha_{j_{1}}|L\cap\overline{U}}=g_{\gamma|L\cap\overline{U}}$ for all
$n\geq j_{0}.$ Thus $\left(  L\cap\overline{U}\right)  ^{1}\left(  \left(
f_{n}\right)  ,\delta\right)  =\emptyset.$ By Lemma \ref{1}(d),%
\[
L^{1}\left(  \left(  f_{n}\right)  ,\delta\right)  \cap\left(  L\cap U\right)
=\emptyset.
\]
In particular,%
\[
x\notin L^{1}\left(  \left(  f_{n}\right)  ,\delta\right)  =H^{\gamma
+1}\left(  \left(  f_{n}\right)  ,\delta\right)  .
\]
Since $x\in H^{\gamma}\left(  \left(  f_{n}\right)  ,\delta\right)  \setminus
K_{\gamma+1}$ is arbitrary, this shows that $H^{\gamma+1}\left(  \left(
f_{n}\right)  ,\delta\right)  \subseteq K_{\gamma+1}.$
\end{proof}

\bigskip

In particular, if $\beta_{K}\left(  f\right)  \leq\beta_{0}\cdot\gamma_{0},$
we have the following.

\begin{theorem}
\label{6}Let $f$ be a Baire-1 function on $K,$ $1\leq\beta_{0},\gamma
_{0}<\omega_{1},$ and $\beta\left(  f\right)  \leq\beta_{0}\cdot\gamma_{0}.$
For any $\varepsilon>0,$ there exist $\tilde{f}:K\rightarrow\mathbb{R}$ and a
sequence of functions $f_{n}:K\rightarrow\mathbb{R}$ such that $\left(
f_{n}\right)  $ converges to $\tilde{f}$ pointwise, $\left\|  \tilde
{f}-f\right\|  \leq\varepsilon$, $\beta\left(  f_{n}\right)  \leq\beta_{0}$
for all $n\in\mathbb{N}$, and $\gamma\left(  \left(  f_{n}\right)  \right)
\leq\gamma_{0}.$
\end{theorem}

A couple more preparatory steps will allow us to improve Theorem \ref{6} to an
exact result (i.e., $\varepsilon=0$) when $\gamma_{0}$ is of the right form.

\begin{theorem}
[{\cite[Lemma 2.5]{K-L}}]\label{9}If $\left(  f_{n}\right)  $ and $\left(
g_{n}\right)  $ are two sequences of real-valued functions on $K$ such that
$\gamma\left(  \left(  f_{n}\right)  \right)  \leq\omega^{\xi}$ and
$\gamma\left(  \left(  g_{n}\right)  \right)  \leq\omega^{\xi}$ for some
$\xi<\omega_{1},$ then $\gamma\left(  \left(  f_{n}+g_{n}\right)  \right)
\leq\omega^{\xi}.$
\end{theorem}

\begin{proposition}
\label{8}For $1\leq\xi<\omega_{1},$ $\frak{B}_{1}^{\xi}\left(  K\right)
=\left\{  f\in\mathbb{R}^{K}:\beta\left(  f\right)  \leq\omega^{\xi}\right\}
$ is a vector subspace of $\mathbb{R}^{K}$ that is closed under the topology
uniform convergence.
\end{proposition}

We postpone the proof of Proposition \ref{8} until the next section. We are
now ready to prove the converse of Theorem \ref{4} in certain cases.

\begin{theorem}
\label{11}If $f\in\frak{B}_{1}\left(  K\right)  $ and $\beta\left(  f\right)
\leq\beta_{0}\cdot\omega^{\gamma_{0}}$ for some $1\leq\beta_{0}<\omega_{1}$
and $\gamma_{0}<\omega_{1},$ then there exists $\left(  f_{n}\right)
\subseteq\frak{B}_{1}\left(  K\right)  $ such that $\left(  f_{n}\right)  $
converges pointwise to $f,$ $\beta\left(  f_{n}\right)  \leq\beta_{0}$ for all
$n\in\mathbb{N}$ and $\gamma\left(  \left(  f_{n}\right)  \right)  \leq
\omega^{\gamma_{0}}.$
\end{theorem}

\begin{proof}
First we assume $\beta_{0}$ is of the form $\omega^{\alpha_{0}},$ where
$\alpha_{0}<\omega_{1}.$ By Theorem \ref{6} there exist a sequence $\left(
f_{n}^{1}\right)  \subseteq\frak{B}_{1}\left(  K\right)  $ and a function
$f^{1}\in\frak{B}_{1}\left(  K\right)  $ such that$,$ $\beta\left(  f_{n}%
^{1}\right)  \leq\omega^{\alpha_{0}}$ for all $n,$ $\left(  f_{n}^{1}\right)
$ converges pointwise to $f^{1},$ $\left\|  f^{1}-f\right\|  \leq\dfrac{1}%
{2},$ and $\gamma\left(  \left(  f_{n}^{1}\right)  \right)  \leq\omega
^{\gamma_{0}}.$ Then $\beta\left(  f^{1}\right)  \leq\omega^{\alpha_{0}}%
\cdot\omega^{\gamma_{0}}=\omega^{\alpha_{0}+\gamma_{0}}$ by Theorem \ref{4}$.$
This implies that $\beta\left(  f-f^{1}\right)  \leq\omega^{\alpha_{0}%
+\gamma_{0}}$ by Proposition \ref{8}$.$ Hence there exist $\left(  f_{n}%
^{2}\right)  \subseteq\frak{B}_{1}\left(  K\right)  $ and $f^{2}$ such that
$\beta\left(  f_{n}^{2}\right)  \leq\omega^{\alpha_{0}}$ for all
$n\in\mathbb{N}$, $\left(  f_{n}^{2}\right)  $ converges pointwise to
$f^{2},\left\|  f-f^{1}-f^{2}\right\|  \leq\dfrac{1}{2^{2}},$ and
$\gamma\left(  \left(  f_{n}^{2}\right)  \right)  \leq\omega^{\gamma_{0}}.$ We
may assume that $\left\|  f_{n}^{2}\right\|  \leq\dfrac{1}{2}$ for all
$n\in\mathbb{N},$ for otherwise, simply replace $f_{n}^{2}$ by $\hat{f}%
_{n}^{2}=\left(  f_{n}^{2}\vee\tfrac{-1}{2}\right)  \wedge\tfrac{1}{2}.$
Continuing, we obtain \ $f^{m}$ and $\left(  f_{n}^{m}\right)  _{n=1}^{\infty
}$ for each $m$ such that

\begin{itemize}
\item $\left\|  f_{n}^{m}\right\|  \leq\dfrac{1}{2^{m-1}},$

\item $\beta\left(  f_{n}^{m}\right)  \leq\omega^{\alpha_{0}}$ for all
$m,\,n\in\mathbb{N},$

\item $\gamma\left(  \left(  f_{n}^{m}\right)  _{n}\right)  \leq\omega
^{\gamma_{0}}$ for all $m\in\mathbb{N},$

\item $f^{m}=\lim\limits_{n}f_{n}^{m}$ (pointwise) for all $m\in\mathbb{N},$ and

\item $\sum_{m=1}^{\infty}f^{m}$ converges uniformly to $f$ on $K.$
\end{itemize}

Let $g_{n}^{m}=f_{n}^{1}+f_{n}^{2}+...+f_{n}^{m}$ and $g_{n}=\sum
_{m=1}^{\infty}f_{n}^{m}.$ By Theorem \ref{9}, $\gamma\left(  \left(
g_{n}^{m}\right)  _{n}\right)  \leq\omega^{\gamma_{0}}$ for all $m\in
\mathbb{N}.$ Given $\varepsilon>0,$ there exists $m_{0}$ such that for all
$n\in\mathbb{N},$ $\left\|  g_{n}^{m_{0}}-g_{n}\right\|  \leq\varepsilon.$
Then $K^{\omega^{\gamma_{0}}}\left(  \left(  g_{n}\right)  ,3\varepsilon
\right)  \subseteq K^{\omega^{\gamma_{0}}}\left(  \left(  g_{n}^{m_{0}%
}\right)  ,\varepsilon\right)  =\emptyset.$ Therefore $\gamma\left(  \left(
g_{n}\right)  \right)  \leq\omega^{\gamma_{0}}.$ By Proposition \ref{8},
$\beta\left(  g_{n}^{m}\right)  \leq\omega^{\alpha_{0}}$ for all $m,\,n.$
Therefore, $\beta\left(  g_{n}\right)  \leq\omega^{\alpha_{0}}$ by Proposition
\ref{8}. Moreover,
\begin{align*}
\lim\limits_{n}g_{n}  &  =\lim\limits_{n}\lim\limits_{m}g_{n}^{m}%
=\lim\limits_{m}\lim\limits_{n}g_{n}^{m}\\
&  =\lim\limits_{m}\sum_{k=1}^{m}f^{k}=f\text{ pointwise.}%
\end{align*}
This proves the theorem in case $\beta_{0}=\omega^{\alpha_{0}}$, with $\left(
g_{n}\right)  $ in place of $\left(  f_{n}\right)  .$

For a general nonzero countable ordinal $\beta_{0},$ write $\beta_{0}$ in
Cantor normal form as
\[
\beta_{0}=\omega^{\beta_{1}}\cdot m_{1}+\omega^{\beta_{2}}\cdot m_{2}%
+...+\omega^{\beta_{k}}\cdot m_{k},
\]
where $k,\,m_{1},...,m_{k}\in\mathbb{N},$ $\omega_{1}>\beta_{1}>\beta
_{2}>...>\beta_{k}.$ If $\gamma_{0}\neq0,$ then $\beta_{0}\cdot\omega
^{\gamma_{0}}=\omega^{\beta_{1}}\cdot\omega^{\gamma_{0}}.$ By the previous
case, there exists $\left(  f_{n}\right)  \subseteq\frak{B}_{1}\left(
K\right)  $ such that $\beta\left(  f_{n}\right)  \leq\omega^{\beta_{1}}%
\leq\beta_{0},$ $\gamma\left(  \left(  f_{n}\right)  \right)  \leq
\omega^{\gamma_{0}}$ and $\left(  f_{n}\right)  $ converges pointwise to $f.$
If $\gamma_{0}=0,$ take $f_{n}=f$ for all $n.$ Then $\beta\left(
f_{n}\right)  \leq\beta_{0}$ for all $n$, $\gamma\left(  \left(  f_{n}\right)
\right)  =1=\omega^{\gamma_{0}}$ and $\left(  f_{n}\right)  $ converges
pointwise to $f$.
\end{proof}

The combination of Theorem \ref{4} and Corollary \ref{11a} yields Theorem
\ref{main}.

\begin{corollary}
\label{11a}Let $f$ $\in\frak{B}_{1}^{\xi}\left(  K\right)  ,$ respectively,
$\mathcal{B}_{1}^{\xi}\left(  K\right)  ,$ for some $\xi<\omega_{1}.$ For all
countable ordinals $\mu,$ $\nu$ such that $\mu+\nu\geq\xi,$ there exists a
sequence $\left(  f_{n}\right)  \subseteq\frak{B}_{1}^{\mu}\left(  K\right)
$, respectively, $\mathcal{B}_{1}^{\mu}\left(  K\right)  ,$ such that
$f_{n}\rightarrow f$ pointwise, and $\gamma\left(  \left(  f_{n}\right)
\right)  \leq\omega^{\nu}.$
\end{corollary}

We do not know if Theorem \ref{11} holds without the restriction on the form
of the ordinal $\gamma\left(  \left(  f_{n}\right)  \right)  .$

\begin{problem}
Is it true that if $f\in\frak{B}_{1}\left(  K\right)  $ with $\beta\left(
f\right)  \leq\beta_{0}\cdot\gamma_{0}$ for some countable ordinals $\beta
_{0}$ and $\gamma_{0},$ then there exists a sequence $\left(  f_{n}\right)  $
converging pointwise to $f$ so that $\sup\limits_{n}\beta\left(  f_{n}\right)
\leq\beta_{0}$ and $\gamma\left(  \left(  f_{n}\right)  \right)  \leq
\gamma_{0}?$
\end{problem}

As another application of our results, we give the proof of another
characterization of the classes $\mathcal{B}_{1}^{\xi}\left(  K\right)  $ due
to Kechris and Louveau.

\begin{definition}
[{\cite[Section 3]{K-L}}]A family $\left\{  \Phi_{\xi}:0\leq\xi<\omega
_{1}\right\}  $ of real-valued functions on $K$ is defined as follows.%
\[
\Phi_{0}=C\left(  K\right)  ,
\]%
\[
\Phi_{\xi+1}=\left\{
\begin{array}
[c]{c}%
f:f\text{ is the pointwise limit of a bounded sequence}\\
\text{ }\left(  f_{n}\right)  \subseteq\Phi_{\xi}\text{ such that }%
\gamma\left(  \left(  f_{n}\right)  \right)  \leq\omega.
\end{array}
\right\}  ,
\]
and for limit ordinals $\lambda,$%
\[
\Phi_{\lambda}=\left\{
\begin{array}
[c]{c}%
f:f\text{ is the uniform limit of a bounded sequence}\\
\text{ }\left(  f_{n}\right)  \subseteq\bigcup_{\xi<\lambda}\Phi_{\xi}\text{.}%
\end{array}
\right\}  .
\]
\end{definition}

\begin{corollary}
[{\cite[Theorem 4.2]{K-L}}]For each $\xi<\omega_{1},$ $\mathcal{B}_{1}^{\xi
}\left(  K\right)  =\Phi_{\xi}.$
\end{corollary}

\begin{proof}
The case $\xi=0$ is trivial. Suppose the corollary holds for some $\xi
<\omega_{1}.$ If $f\in\mathcal{B}_{1}^{\xi+1}\left(  K\right)  ,$ it follows
from Corollary \ref{11a} that $f$ is the pointwise limit of a bounded sequence
$\left(  f_{n}\right)  $ in $\mathcal{B}_{1}^{\xi}\left(  K\right)  $ such
that $\gamma\left(  \left(  f_{n}\right)  \right)  \leq\omega.$ Since
$\mathcal{B}_{1}^{\xi}\left(  K\right)  =\Phi_{\xi}$ by the inductive
hypothesis, $f\in\Phi_{\xi+1}.$ Conversely, if $f\in\Phi_{\xi+1},$ then $f$ is
the pointwise limit of a sequence $\left(  f_{n}\right)  $ in $\Phi_{\xi}$
with $\gamma\left(  \left(  f_{n}\right)  \right)  \leq\omega.$ Since
$\Phi_{\xi}=\mathcal{B}_{1}^{\xi}\left(  K\right)  ,$ $\beta\left(  f\right)
\leq\omega^{\xi+1}$ by Theorem \ref{4}$.$ Thus $f\in\mathcal{B}_{1}^{\xi
+1}\left(  K\right)  .$

Now assume that the corollary holds for all $\xi^{\prime}<\xi,$ where $\xi$ is
a countable limit ordinal. Let $f\in\Phi_{\xi}.$ By the inductive hypothesis,
$\Phi_{\xi^{\prime}}=\mathcal{B}_{1}^{\xi^{\prime}}\left(  K\right)
\subseteq\mathcal{B}_{1}^{\xi}\left(  K\right)  $ for $\xi^{\prime}<\xi.$
Hence $f$ is the uniform limit of a sequence in $\mathcal{B}_{1}^{\xi}\left(
K\right)  $, and thus belongs to $\mathcal{B}_{1}^{\xi}\left(  K\right)  .$
Conversely, assume that $f\in\mathcal{B}_{1}^{\xi}\left(  K\right)  .$ For
every $n\in\mathbb{N},$ there exists $\xi_{n}<\xi$ such that $\beta\left(
f,\frac{1}{n}\right)  \leq\omega^{\xi_{n}}.$ By Corollary \ref{c9}, the exists
$f_{n}\in\mathcal{B}_{1}^{\xi_{n}}\left(  K\right)  =\Phi_{\xi_{n}}$ such that
$\left\|  f-f_{n}\right\|  \leq\frac{1}{n}.$ Thus $f\in\Phi_{\xi},$ as required.
\end{proof}

\begin{remark}
\emph{If a family }$\left\{  \Psi_{\xi}:0\leq\xi<\omega_{1}\right\}  $\emph{
is defined in a similar way as the family }$\left\{  \Phi_{\xi}:0\leq
\xi<\omega_{1}\right\}  $\emph{ except for the removal of the boundedness
condition on the sequence }$\left(  f_{n}\right)  ,$\emph{ then }$\Psi_{\xi
}=B_{1}^{\xi}\left(  K\right)  $\emph{ for all }$\xi<\omega_{1}.$
\end{remark}

\section{Optimal limit of continuous functions}

In this section we prove the equivalence of the indices $\beta$ and $\gamma$
for functions in $\frak{B}_{1}\left(  K\right)  $ in the same sense that was
established for $\mathcal{B}_{1}\left(  K\right)  $ in Theorem 2.3 of
\cite{K-L}. Namely, it is shown that for all $f\in\frak{B}_{1}\left(
K\right)  ,$ $\beta\left(  f\right)  $ is the smallest ordinal $\gamma_{0}$
for which there exists a sequence $\left(  f_{n}\right)  $ in $C\left(
K\right)  $ converging pointwise to $f$ and satisfying $\gamma\left(  \left(
f_{n}\right)  \right)  \leq\gamma_{0}.$ Let us note that this result is also
the converse of Theorem \ref{4} when $\beta_{0}=1.$

\begin{definition}
Let $\left(  f_{n}\right)  \subseteq\mathbb{R}^{K}$ and $f\in\mathbb{R}^{K}.$
We write

(a) $\left(  g_{n}\right)  \prec\left(  f_{n}\right)  $ if $\left(
g_{n}\right)  $ is a convex block combination of $\left(  f_{n}\right)  ,$
i.e., there exists a sequence of non-negative real numbers $\left(
a_{k}\right)  $ and a strictly increasing sequence $\left(  p_{n}\right)  $ in
$\mathbb{N}$ such that $\sum_{k=p_{n-1}+1}^{p_{n}}a_{k}=1$ and $g_{n}%
=\sum_{k=p_{n-1}+1}^{p_{n}}a_{k}f_{k}$ for all $n$ $\left(  p_{0}=0\right)  .$

(b) $\left(  g_{n}\right)  \overset{a}{\prec}\left(  f_{n}\right)  $ if there
exists $m\in\mathbb{N}$ such that $\left(  g_{n}\right)  _{n=m}^{\infty}%
\prec\left(  f_{n}\right)  ,$ and

(c) $\left[  f\right]  _{-M}^{M}=\left(  f\vee-M\right)  \wedge M,$ where
$0\leq M\in\mathbb{R}.$
\end{definition}

The easy proof of the next lemma is left to the reader.

\begin{lemma}
\label{14}If $\left(  g_{n}\right)  \overset{a}{\prec}\left(  f_{n}\right)  ,$
then $\gamma\left(  \left(  g_{n}\right)  ,\varepsilon\right)  \leq
\gamma\left(  \left(  f_{n}\right)  ,\varepsilon\right)  $ for all
$\varepsilon>0.$
\end{lemma}

\begin{lemma}
\label{15}Let $f$ be a Baire-1 function on $K.$ Suppose $\mathcal{H}$ is a
countable collection of compact subsets of $K$ such that $\left\|  f\right\|
_{H}<\infty$ for all $H\in\mathcal{H}$ and $\bigcup_{H\in\mathcal{H}}H=K.$
Then there exists $\left(  f_{n}\right)  \subseteq C\left(  K\right)  $ such that

(i) $f_{n}$ $\rightarrow$ $f$ pointwise, and

(ii) $\left(  f_{n|H}\right)  $ is a bounded subset of $C\left(  H\right)  $
for all $H\in\mathcal{H}.$
\end{lemma}

\begin{proof}
Write $\mathcal{H}$ as a sequence $\left(  H_{m}\right)  _{m=1}^{\infty}.$
Without loss of generality, assume that $H_{m}\subseteq H_{m+1}$ for all
$m\in\mathbb{N}.$ Since $f$ is Baire-1, there exists $\left(  f_{n}%
^{0}\right)  \subseteq C\left(  K\right)  $ such that $\left(  f_{n}%
^{0}\right)  $ converges pointwise to $f.$ Assume that $\left(  f_{n}%
^{m-1}\right)  _{n}\subseteq C\left(  K\right)  $ has been chosen so that
$\lim\limits_{n}f_{n}^{m-1}=f$ pointwise. If $m,\,n\in\mathbb{N}$, let
$U_{n}^{m}$ be the $\dfrac{1}{n}-$neighborhood of $H_{m}$ in $K$ and let
$M_{m}=\left\|  f\right\|  _{H_{m}}.$ For all $n,$ the function $\left[
f_{n}^{m-1}\right]  _{-M_{m}|H_{m}}^{M_{m}}\cup f_{n|K\setminus U_{n}^{m}%
}^{m-1}$ is continuous on $H_{m}\cup\left(  K\setminus U_{n}^{m}\right)  .$
Let $f_{n}^{m}$ be a continuous extension of the function onto $K.$ Then
$\left(  f_{n}^{m}\right)  \subseteq C\left(  K\right)  .$ If $x\in H_{m},$
then $\lim\limits_{n}f_{n}^{m}\left(  x\right)  =\lim\limits_{n}\left[
f_{n}^{m-1}\left(  x\right)  \right]  _{-M_{m}}^{M_{m}}=\left[  f\left(
x\right)  \right]  _{-M_{m}}^{M_{m}}=f\left(  x\right)  $ since $\left\|
f\right\|  _{H_{m}}=M_{m}.$ If $x\notin H_{m},$ then there exists $n_{0}$ such
that $x\in K\setminus U_{n_{0}}^{m};$ thus $x\in K\setminus U_{n}^{m}$ for all
$n\geq n_{0}.$ Therefore $f_{n}^{m}\left(  x\right)  =f_{n}^{m-1}\left(
x\right)  $ for all $n\geq n_{0}.$ Hence $\lim\limits_{n}f_{n}^{m}\left(
x\right)  =f\left(  x\right)  .$ Thus $\lim\limits_{n}f_{n}^{m}=f$
pointwise$.$ Now for each $n\in\mathbb{N},$ let $f_{n}=f_{n}^{n}.$ Since
$H_{m}\subseteq H_{n}$ for all $n\geq m$, on $H_{m}$ we have
\begin{align*}
f_{n}  &  =f_{n}^{n}=\left[  f_{n}^{n-1}\right]  _{-M_{n}}^{M_{n}}\\
&  =\left[  \left[  f_{n}^{n-2}\right]  _{-M_{n-1}}^{M_{n-1}}\right]
_{-M_{n}}^{M_{n}}=...=\left[  ...\left[  \left[  f_{n}^{m-1}\right]  _{-M_{m}%
}^{M_{m}}\right]  _{-M_{m+1}}^{M_{m+1}}...\right]  _{-M_{n}}^{M_{n}}\text{ }\\
&  =\left[  f_{n}^{m-1}\right]  _{-M_{m}}^{M_{m}}\text{ as }M_{m}\leq
M_{m+1}\leq...\leq M_{n}.
\end{align*}
Thus $f_{n}=\left[  f_{n}^{m-1}\right]  _{-M_{m}}^{M_{m}}$ on $H_{m}$ for all
$n\geq m.$ In particular, on the set $H_{m},$
\[
\lim\limits_{n}f_{n}=\left[  \lim\limits_{n}f_{n}^{m-1}\right]  _{-M_{m}%
}^{M_{m}}=\left[  f\right]  _{-M_{m}}^{M_{m}}=f
\]
since $\left\|  f\right\|  _{H_{m}}=M_{m}.$ As $K=\bigcup H_{m},$ we see that
$f_{n}\rightarrow f$ pointwise$.$ Also, for each $m,$ $\left(  f_{n|H_{m}%
}\right)  _{n=m}^{\infty}$ is bounded (by $M_{m}$) in $C\left(  H_{m}\right)
;$ thus $\left(  f_{n|H_{m}}\right)  _{n=1}^{\infty}$ is bounded in $C\left(
H_{m}\right)  .$
\end{proof}

For the next lemma, recall that for a real-valued function $f$ defined on a
set $S,$ $\operatorname*{osc}\left(  f,S\right)  =\sup\left\{  \left|
f\left(  s_{1}\right)  -f\left(  s_{2}\right)  \right|  :s_{1},\,s_{2}\in
S\right\}  .$

\begin{lemma}
\label{16}Let $\left(  f_{n}\right)  $ be bounded in $C\left(  H\right)  ,$
where $H$ is a compact metric space. Suppose $\left(  f_{n}\right)  $
converges pointwise to $f$ and $H^{1}\left(  f,\varepsilon\right)  =\emptyset$
for some $\varepsilon>0,$ then there exists $\left(  g_{n}\right)
\prec\left(  f_{n}\right)  $ such that $H^{1}\left(  \left(  g_{n}\right)
,7\varepsilon\right)  =\emptyset.$
\end{lemma}

\begin{proof}
By Corollary \ref{c9}, there exists $\tilde{f}\in C\left(  H\right)  $ such
that $\left\|  f-\tilde{f}\right\|  _{H}\leq\varepsilon.$ Then $\left(
f_{n}-\tilde{f}\right)  $ is bounded in $C\left(  H\right)  ,$ $f_{n}%
-\tilde{f}\rightarrow f-\tilde{f}$ pointwise and $\operatorname*{osc}\left(
f-\tilde{f},H\right)  \leq2\varepsilon$. By the first statement in the proof
of Theorem 2.3 in \cite{K-L}, there exists $\left(  h_{n}\right)  \prec\left(
f_{n}-\tilde{f}\right)  $ such that $\left\|  h_{n}-(f-\tilde{f})\right\|
_{H}\leq3\varepsilon.$ Let $g_{n}=h_{n}+\tilde{f}$ for all $n\in\mathbb{N}$.
Then $\left(  g_{n}\right)  \prec\left(  f_{n}\right)  $ and $\left\|
g_{n}-f\right\|  _{H}\leq3\varepsilon$ for all $n\in\mathbb{N}$. It follows
that $H^{1}\left(  \left(  g_{n}\right)  ,7\varepsilon\right)  =\emptyset.$
\end{proof}

\begin{theorem}
\label{18}Let $f$ be a Baire-1 function on $K.$ There exists a sequence
$\left(  f_{n}\right)  \subseteq C\left(  K\right)  $ such that $\left(
f_{n}\right)  $ converges pointwise to $f$ and $\gamma\left(  \left(
f_{n}\right)  \right)  =\beta\left(  f\right)  .$
\end{theorem}

\begin{proof}
Let $\beta_{0}=\beta\left(  f\right)  .$ For each $\alpha<\beta_{0},$ and all
$m,\,j\in\mathbb{N}$, let $U_{m,j}^{\alpha}$ be the $\dfrac{1}{j}%
-$neighborhood of $K^{\alpha}\left(  f,\dfrac{1}{m}\right)  $ in $K.$ Define
\[
\mathcal{H}=\left\{  K^{\alpha}\left(  f,\dfrac{1}{m}\right)  \setminus
U_{m,j}^{\alpha+1}:\alpha<\beta_{0},\,m,\,j\in\mathbb{N}\right\}  .
\]
Then $\mathcal{H}$ is a countable collection of compact subsets of $K$ such
that $\bigcup_{H\in\mathcal{H}}H=K.$ If $\alpha<\beta_{0}$ and $m,\,j\in
\mathbb{N}$, by Lemma \ref{5}, there is a continuous function $g$ on
$H=K^{\alpha}\left(  f,\dfrac{1}{m}\right)  \setminus U_{m,j}^{\alpha+1}$ such
that $\left\|  g-f\right\|  _{H}\leq\dfrac{1}{m}.$ Hence $\left\|  f\right\|
_{H}<\infty$ for all $H\in\mathcal{H}.$ By Lemma \ref{15}, there exists
$\left(  g_{n}\right)  \subseteq C\left(  K\right)  $ such that $\left(
g_{n}\right)  $ converges pointwise to $f$ and $\left(  g_{n|H}\right)  $ is
bounded in $C\left(  H\right)  $ for all $H\in\mathcal{H}$.

List the elements of $\mathcal{H}$ in a sequence $\left(  H_{k}\right)
_{k=1}^{\infty}.$ Take $\varepsilon_{k}=\dfrac{1}{m}$ if $H_{k}$ is of the
form $K^{\alpha}\left(  f,\dfrac{1}{m}\right)  \setminus U_{m,j}^{\alpha+1}$
for some $\alpha,\,m,\,j.$ Let $\left(  g_{n}^{0}\right)  =\left(
g_{n}\right)  .$ Suppose $\left(  g_{n}^{k-1}\right)  _{n}\prec\left(
g_{n}\right)  _{n}$ has been chosen. Then $\left(  g_{n}^{k-1}\right)  _{n}$
converges to $f$ pointwise, $\left(  g_{n|H_{k}}^{k-1}\right)  $ is a bounded
sequence in $C\left(  H_{k}\right)  ,$ and $\left(  H_{k}\right)  ^{1}\left(
f,\varepsilon_{k}\right)  =\emptyset.$ By Lemma \ref{16}, there exists
$\left(  g_{n}^{k}\right)  _{n}\prec\left(  g_{n}^{k-1}\right)  _{n}$ such
that $\left(  H_{k}\right)  ^{1}\left(  \left(  g_{n}^{k}\right)
_{n},7\varepsilon_{k}\right)  =\emptyset.$ Let $f_{n}=g_{n}^{n}$ for all
$n\in\mathbb{N}$. Then $\left(  f_{n}\right)  \prec\left(  g_{n}\right)  .$
Therefore $\left(  f_{n}\right)  \subseteq C\left(  K\right)  $ and $\left(
f_{n}\right)  $ converges pointwise to $f$. We claim that for all
$m\in\mathbb{N}$ and for all $\alpha\leq\beta_{0},$ $K^{\alpha}\left(  \left(
f_{n}\right)  ,\dfrac{7}{m}\right)  \subseteq K^{\alpha}\left(  f,\dfrac{1}%
{m}\right)  .$ We prove the claim by induction on $\alpha.$ The claim is
trivial if $\alpha=0$ or $\alpha$ is a limit ordinal. Assume that $\alpha
\leq\beta_{0}$ is a successor ordinal and that the claim holds for $\alpha-1.$
Let $x\in K^{\alpha}\left(  \left(  f_{n}\right)  ,\dfrac{7}{m}\right)  .$
Then $x\in K^{\alpha-1}\left(  \left(  f_{n}\right)  ,\dfrac{7}{m}\right)
\subseteq K^{\alpha-1}\left(  f,\dfrac{1}{m}\right)  .$ If $x\notin K^{\alpha
}\left(  f,\dfrac{1}{m}\right)  ,$ there exists $j\in\mathbb{N}$ such that
$d\left(  x,K^{\alpha}\left(  f,\dfrac{1}{m}\right)  \right)  >\dfrac{1}{j}.$
Choose $k$ such that $H_{k}=K^{\alpha-1}\left(  f,\dfrac{1}{m}\right)
\setminus U_{m,j}^{\alpha}.$ Then $\left(  f_{n}\right)  \overset{a}{\prec
}\left(  g_{n}^{k}\right)  _{n}$ and $\gamma_{H_{k}}\left(  \left(  g_{n}%
^{k}\right)  _{n},7\varepsilon_{k}\right)  \leq1$ since $\left(  H_{k}\right)
^{1}\left(  \left(  g_{n}^{k}\right)  _{n},7\varepsilon_{k}\right)
=\emptyset.$ By Lemma \ref{14}, $\left(  H_{k}\right)  ^{1}\left(  \left(
f_{n}\right)  ,7\varepsilon_{k}\right)  =\emptyset.$ Thus $\left(
H_{k}\right)  ^{1}\left(  \left(  f_{n}\right)  ,\dfrac{7}{m}\right)
=\emptyset.$ But since $d\left(  x,K^{\alpha}\left(  f,\dfrac{1}{m}\right)
\right)  >\dfrac{1}{j},$ there exists an open set $U$ in $\tilde{K}%
=K^{\alpha-1}\left(  f,\dfrac{1}{m}\right)  $ such that $x\in U\subseteq
H_{k}\subseteq\tilde{K}.$ By Lemma \ref{1}(d), $\left(  \tilde{K}\right)
^{1}\left(  \left(  f_{n}\right)  ,\dfrac{7}{m}\right)  \cap U\subseteq\left(
H_{k}\right)  ^{1}\left(  \left(  f_{n}\right)  ,\dfrac{7}{m}\right)
=\emptyset.$ Therefore $x\notin\left(  \tilde{K}\right)  ^{1}\left(  \left(
f_{n}\right)  ,\dfrac{7}{m}\right)  =K^{\alpha}\left(  \left(  f_{n}\right)
,\dfrac{7}{m}\right)  ,$ a contradiction. This proves the claim. From the
claim, $K^{\beta_{0}}\left(  \left(  f_{n}\right)  ,\dfrac{7}{m}\right)
\subseteq K^{\beta_{0}}\left(  f,\dfrac{1}{m}\right)  =\emptyset$ for all
$m\in\mathbb{N}.$ Therefore $\gamma\left(  \left(  f_{n}\right)  \right)
\leq\beta_{0}.$ Since $\gamma\left(  \left(  f_{n}\right)  \right)  \geq
\beta_{0}$ by \cite[Proposition 2.1]{K-L}, (or Theorem \ref{4}),
$\gamma\left(  \left(  f_{n}\right)  \right)  =\beta_{0}=\beta\left(
f\right)  .$
\end{proof}

\begin{remark}
\emph{Unlike in Theorem 2.3 of \cite{K-L}, in general we cannot get a sequence
}$\left(  g_{n}\right)  \prec\left(  f_{n}\right)  $\emph{ such that }%
$\gamma\left(  \left(  g_{n}\right)  \right)  =\beta\left(  f\right)  .$\emph{
Indeed, let }$K=\left[  0,1\right]  $\emph{ and for each }$n\in N$\emph{ let
}$f_{n}$\emph{ be a continuous function that vanishes outside }$\left[
\dfrac{1}{n+1},\dfrac{1}{n}\right]  $\emph{ such that }$\int_{K}f_{n}%
=1.$\emph{ Then }$\left(  f_{n}\right)  $\emph{ converges pointwise to }%
$f=0$\emph{. Suppose }$\left(  g_{n}\right)  \prec\left(  f_{n}\right)
,$\emph{ then }$\int_{K}g_{n}=1$\emph{ for all }$n\in N.$\emph{ Thus }$\left(
g_{n}\right)  $\emph{ does not converge uniformly to }$f,$\emph{ i.e.,
}$\gamma\left(  \left(  g_{n}\right)  \right)  >1=\beta\left(  f\right)  .$
\end{remark}

\begin{proof}
[Proof of Proposition \ref{8}]It is easy to see that for all $f\in\frak{B}%
_{1}^{\xi}\left(  K\right)  $ and $a\in\mathbb{R},$ $af\in\frak{B}_{1}^{\xi
}\left(  K\right)  .$ If $f$,\thinspace$g\in\frak{B}_{1}^{\xi}\left(
K\right)  ,$ then by Theorem \ref{18} there exist two sequences of continuous
functions $\left(  f_{n}\right)  $ and $\left(  g_{n}\right)  $ converging
pointwise to $f$ and $g$ respectively such that $\gamma\left(  \left(
f_{n}\right)  \right)  \leq\omega^{\xi}$ and $\gamma\left(  \left(
g_{n}\right)  \right)  \leq\omega^{\xi}.$ According to Theorem \ref{9},
$\gamma\left(  \left(  f_{n}+g_{n}\right)  \right)  \leq\omega^{\xi}.$ Hence
by Theorem \ref{4}, $f+g\in\frak{B}_{1}^{\xi}\left(  K\right)  .$ Finally,
given $f\in\overline{\frak{B}_{1}^{\xi}\left(  K\right)  }$ and $\varepsilon
>0$, choose $g\in\frak{B}_{1}^{\xi}\left(  K\right)  $ such that $\left\|
f-g\right\|  \leq\dfrac{\varepsilon}{3}.$ Then $K^{\omega^{\xi}}\left(
f,\varepsilon\right)  \subseteq K^{\omega^{\xi}}\left(  g,\dfrac{\varepsilon
}{3}\right)  =\emptyset.$ Thus $f\in\frak{B}_{1}^{\xi}\left(  K\right)  .$
\end{proof}

\section{Product of Baire-1 functions}

In \cite{K-L}, it is observed that the classes $\mathcal{B}_{1}^{\xi}\left(
K\right)  ,$ $\xi<\omega_{1}$ are closed under multiplication. However, it is
relative easy to see that this fails for the classes $\frak{B}_{1}^{\xi
}\left(  K\right)  .$ In this section, we show that if $f\in\frak{B}_{1}%
^{\xi_{1}}\left(  K\right)  $ and $g\in\frak{B}_{1}^{\xi_{2}}\left(  K\right)
,$ then $fg\in\frak{B}_{1}^{\xi}\left(  K\right)  ,$ where $\xi=\max\left\{
\xi_{1}+\xi_{2},\xi_{2}+\xi_{1}\right\}  .$ It is also shown that the result
is sharp. The proof of the next lemma is left to the reader.

\begin{lemma}
\label{P1}If $f$ is bounded and $\gamma\left(  \left(  g_{n}\right)  \right)
\leq\xi,$ then $\gamma\left(  \left(  fg_{n}\right)  \right)  \leq\xi.$
\end{lemma}

\begin{lemma}
\label{P2}If $f\in\mathcal{B}_{1}^{\xi_{1}}\left(  K\right)  $ and
$g\in\frak{B}_{1}^{\xi_{2}}\left(  K\right)  ,$ then $fg\in\frak{B}_{1}%
^{\xi_{1}+\xi_{2}}\left(  K\right)  .$
\end{lemma}

\begin{proof}
By Theorem \ref{18}, there exists a sequence $\left(  g_{n}\right)  \subseteq
C\left(  K\right)  $ converging to $g$ pointwise such that $\gamma\left(
\left(  g_{n}\right)  \right)  =\omega^{\xi_{2}}.$ For each $n\in\mathbb{N},$
$g_{n}\in C\left(  K\right)  \subseteq\mathcal{B}_{1}^{\xi_{1}}\left(
K\right)  $ and $f\in\mathcal{B}_{1}^{\xi_{1}}\left(  K\right)  .$ By
\cite{K-L} (see the remark on \cite[p. 217]{K-L}), $fg_{n}\in\mathcal{B}%
_{1}^{\xi_{1}}\left(  K\right)  .$ Lemma \ref{P1} implies that $\gamma\left(
\left(  fg_{n}\right)  \right)  \leq\omega^{\xi_{2}}.$ Since $\left(
fg_{n}\right)  $ converges to $fg$ pointwise, it follows from Theorem \ref{4}
that $\beta\left(  fg\right)  \leq\omega^{\xi_{1}+\xi_{2}},$ i.e.,
$fg\in\frak{B}_{1}^{\xi_{1}+\xi_{2}}\left(  K\right)  .$
\end{proof}

Now suppose $f\in\frak{B}_{1}^{\xi_{1}}\left(  K\right)  $ and $g\in
\frak{B}_{1}^{\xi_{2}}\left(  K\right)  .$ By Lemma \ref{5}, for all
$\alpha<\omega^{\xi_{2}}$, there is a continuous function $g_{\alpha
}:K^{\alpha}\left(  g,1\right)  \setminus K^{\alpha+1}\left(  g,1\right)
\rightarrow\mathbb{R}$ such that
\[
\left\|  g_{\alpha}-g\right\|  _{K^{\alpha}\left(  g,1\right)  \setminus
K^{\alpha+1}\left(  g,1\right)  }\leq1.
\]
Let $h=\bigcup_{\alpha<\omega^{\xi_{2}}}g_{\alpha}.$ It follows from the proof
of Theorem \ref{AE} that $\beta\left(  h\right)  \leq\omega^{\xi_{2}}.$ Given
a closed set $H\subseteq K,$ we write
\[
\operatorname*{d}\nolimits_{f}\left(  H\right)  =\left\{  x\in H:\underset
{y\in H}{\limsup\limits_{y\rightarrow x}}\left|  f\left(  y\right)  \right|
=\infty\right\}  .
\]
It is easy to see that $\operatorname*{d}\nolimits_{f}\left(  H\right)  $ is a
closed subset of $H$ such that $\operatorname*{d}\nolimits_{f}\left(
H\right)  \subseteq H^{1}\left(  f,\varepsilon\right)  $ for any $\varepsilon>0.$

\begin{lemma}
\label{P3}Suppose that $\alpha<\omega_{1},$ $\delta>0$ and $s>2.$ If
$x\in\left[  K\setminus K^{1}\left(  g,1\right)  \right]  \cap K^{\alpha
}\left(  fh,\delta\right)  ,$ then $x\in K^{\alpha}\left(  f,\dfrac{\delta
}{s\left(  \left|  h\left(  x\right)  \right|  +1\right)  }\wedge1\right)  .$
\end{lemma}

\begin{proof}
The proof is by induction on $\alpha.$ The result is clear if $\alpha=0$ or a
limit ordinal. Assume that the lemma holds for some $\alpha<\omega_{1}.$
Suppose $\delta>0$ and $s>2$ are given. Let $x\in\left[  K\setminus
K^{1}\left(  g,1\right)  \right]  \cap K^{\alpha+1}\left(  fh,\delta\right)
.$ If $x\in\operatorname*{d}_{f}\left(  K^{\alpha}\left(  f,\dfrac{\delta
}{s\left(  \left|  h\left(  x\right)  \right|  +1\right)  }\wedge1\right)
\right)  ,$ then $x\in K^{\alpha+1}\left(  f,\dfrac{\delta}{s\left(  \left|
h\left(  x\right)  \right|  +1\right)  }\wedge1\right)  $ and we are done.
Otherwise, assume that $x\not \in\operatorname*{d}\nolimits_{f}\left(
K^{\alpha}\left(  f,\dfrac{\delta}{s\left(  \left|  h\left(  x\right)
\right|  +1\right)  }\wedge1\right)  \right)  .$ Then there exist a
neighborhood $U_{1}$ of $x$ in $K$ and $M<\infty$ such that $\left|  f\left(
y\right)  \right|  \leq M$ for all $y\in U_{1}\cap K^{\alpha}\left(
f,\dfrac{\delta}{s\left(  \left|  h\left(  x\right)  \right|  +1\right)
}\wedge1\right)  .$ Since $h=g_{0}$ on $K\setminus K^{1}\left(  g,1\right)  ,$
and $g_{0}$ is continuous on $K\setminus K^{1}\left(  g,1\right)  ,$ there
exists a neighborhood $U_{2}$ of $x$ in $K$ such that $\left|  h\left(
x_{1}\right)  -h\left(  x_{2}\right)  \right|  \leq\dfrac{\delta}{2M}$ and
$2\left(  \left|  h\left(  x_{1}\right)  \right|  +1\right)  <s\left(  \left|
h\left(  x\right)  \right|  +1\right)  $ for all $x_{1},x_{2}\in U_{2}.$ Set
$U=\left(  U_{1}\cap U_{2}\right)  \setminus K^{1}\left(  g,1\right)  .$ Then
$U$ is a neighborhood of $x.$

\noindent\textbf{Claim. }$K^{\alpha}\left(  fh,\delta\right)  \cap U\subseteq
K^{\alpha}\left(  f,\dfrac{\delta}{s\left(  \left|  h\left(  x\right)
\right|  +1\right)  }\wedge1\right)  .$

\noindent Note that if $y\in U,$ then $y\in U_{2}.$ Hence there exists $t>2$
such that $t\left(  \left|  h\left(  y\right)  \right|  +1\right)  \leq
s\left(  \left|  h\left(  x\right)  \right|  +1\right)  .$ Also, $y\in
K^{\alpha}\left(  fh,\delta\right)  \cap U$ implies that $y\in\left[
K\setminus K^{1}\left(  g,1\right)  \right]  \cap K^{\alpha}\left(
fh,\delta\right)  .$ Thus $y\in K^{\alpha}\left(  f,\dfrac{\delta}{t\left(
\left|  h\left(  y\right)  \right|  +1\right)  }\wedge1\right)  $ by the
inductive hypothesis. Since
\[
\dfrac{\delta}{t\left(  \left|  h\left(  y\right)  \right|  +1\right)  }%
\geq\dfrac{\delta}{s\left(  \left|  h\left(  x\right)  \right|  +1\right)
}\wedge1,
\]
$y\in K^{\alpha}\left(  f,\dfrac{\delta}{s\left(  \left|  h\left(  x\right)
\right|  +1\right)  }\wedge1\right)  ,$ as required.

Now if $V$ is a neighborhood of $x$ in $K,$ there exist $x_{1},x_{2}\in U\cap
V\cap$\textbf{ }$K^{\alpha}\left(  fh,\delta\right)  $ such that
\begin{align*}
\delta &  \leq\left|  f\left(  x_{1}\right)  h\left(  x_{1}\right)  -f\left(
x_{2}\right)  h\left(  x_{2}\right)  \right| \\
&  \leq\left|  f\left(  x_{1}\right)  -f\left(  x_{2}\right)  \right|  \left|
h\left(  x_{1}\right)  \right|  +\left|  h\left(  x_{1}\right)  -h\left(
x_{2}\right)  \right|  \left|  f\left(  x_{2}\right)  \right| \\
&  \leq\left|  f\left(  x_{1}\right)  -f\left(  x_{2}\right)  \right|  \left|
h\left(  x_{1}\right)  \right|  +\dfrac{\delta}{2M}\cdot M,
\end{align*}
where, in the last inequality, $\left|  f\left(  x_{2}\right)  \right|  \leq
M$ since $x_{2}\in U\cap K^{\alpha}\left(  f,\dfrac{\delta}{s\left(  \left|
h\left(  x\right)  \right|  +1\right)  }\wedge1\right)  $ by the claim.
Therefore,
\[
\left|  f\left(  x_{1}\right)  -f\left(  x_{2}\right)  \right|  \geq
\dfrac{\delta}{s\left(  \left|  h\left(  x\right)  \right|  +1\right)  }%
\wedge1.
\]
By the claim, $x_{1},x_{2}\in V\cap K^{\alpha}\left(  f,\dfrac{\delta
}{s\left(  \left|  h\left(  x\right)  \right|  +1\right)  }\wedge1\right)  .$
Since $V$ is an arbitrary neighborhood of $x,$ this shows that
\[
x\in K^{\alpha+1}\left(  f,\dfrac{\delta}{s\left(  \left|  h\left(  x\right)
\right|  +1\right)  }\wedge1\right)  .
\]
This completes the induction.
\end{proof}

It follows from Lemma \ref{P3} that
\[
K^{\omega^{\xi_{1}}}\left(  fh,\delta\right)  \subseteq K^{1}\left(
g,1\right)  .
\]
Repeating the argument in Lemma \ref{P3} inductively yields

\begin{lemma}
$K^{\omega^{\xi_{1}}\cdot\alpha}\left(  fh,\delta\right)  \subseteq K^{\alpha
}\left(  g,1\right)  $ for all $\alpha<\omega_{1}$, and $\delta>0.$
\end{lemma}

In particular, $K^{\omega^{\xi_{1}}\cdot\omega^{\xi_{2}}}\left(
fh,\delta\right)  =\emptyset$ for all $\delta>0,$ i.e., $fh\in\frak{B}%
_{1}^{\xi_{1}+\xi_{2}}\left(  K\right)  .$

\begin{theorem}
\label{Product}If $f\in\frak{B}_{1}^{\xi_{1}}\left(  K\right)  $ and
$g\in\frak{B}_{1}^{\xi_{2}}\left(  K\right)  ,$ then $fg\in\frak{B}_{1}^{\xi
}\left(  K\right)  ,$ where $\xi=\max\left\{  \xi_{1}+\xi_{2},\,\xi_{2}%
+\xi_{1}\right\}  .$
\end{theorem}

\begin{proof}
From the above, we obtain a function $h$ in $K$ such that $\left\|
g-h\right\|  \leq1,$ $\beta\left(  h\right)  \leq\omega^{\xi_{2}}$ and
$fh\in\frak{B}_{1}^{\xi_{1}+\xi_{2}}\left(  K\right)  .$ Since $g,h\in
\frak{B}_{1}^{\xi_{2}}\left(  K\right)  ,$ it follows from Proposition \ref{8}
that $g-h\in$ $\frak{B}_{1}^{\xi_{2}}\left(  K\right)  .$ As $g-h$ is bounded,
we see that $g-h\in\mathcal{B}_{1}^{\xi_{2}}\left(  K\right)  .$ By Lemma
\ref{P2}, $\left(  g-h\right)  f\in\frak{B}_{1}^{\xi_{2}+\xi_{1}}\left(
K\right)  \subseteq\frak{B}_{1}^{\xi}\left(  K\right)  .$ Also, $fh\in
\frak{B}_{1}^{\xi_{1}+\xi_{2}}\left(  K\right)  \subseteq\frak{B}_{1}^{\xi
}\left(  K\right)  .$ Applying Proposition \ref{8} again gives $fg=f\left(
g-h\right)  +fh\in\frak{B}_{1}^{\xi}\left(  K\right)  .$
\end{proof}

Our final result shows that Theorem \ref{Product} is sharp. We omit the easy
proof of the next lemma.

\begin{lemma}
\label{P4}Suppose that $h\in\frak{B}_{1}\left(  K\right)  ,$ $\alpha
<\omega_{1},$ and $\varepsilon>0.$ Let $V=K\setminus K^{\alpha}\left(
h,\varepsilon\right)  .$ For any $\eta<\omega_{1},$%
\[
K^{\eta}\left(  h,\varepsilon\right)  \setminus K^{\alpha}\left(
h,\varepsilon\right)  \subseteq K^{\eta}\left(  h\chi_{V},\varepsilon\right)
.
\]
\end{lemma}

\begin{theorem}
Suppose that $\xi_{1}$, $\xi_{2}$ are countable ordinals, and let
\[
\xi=\max\left\{  \xi_{1}+\xi_{2},\,\xi_{2}+\xi_{1}\right\}  .
\]
If $K$ is a compact metric space such that $K^{(  \omega^{\xi})
}\neq\emptyset,$ then
\[
\sup\left\{  \beta\left(  fg\right)  :f\in\frak{B}_{1}^{\xi_{1}}\left(
K\right)  ,\,g\in\frak{B}_{1}^{\xi_{2}}\left(  K\right)  \right\}
=\omega^{\xi}.
\]
\end{theorem}

\begin{proof}
We may of course assume that neither $\xi_{1}$ nor $\xi_{2}$ is $0,$ and that
$\xi=\xi_{1}+\xi_{2}.$ The assumption on $K$ yields a $\left\{  0,1\right\}
$-valued function $h$ in $\frak{B}_{1}\left(  K\right)  $ such that
$K^{\omega^{\xi}}\left(  h,1\right)  \neq\emptyset.$ Denote $K^{\alpha}\left(
h,1\right)  $ by $K_{\alpha},$ $\alpha<\omega_{1}.$ Choose a sequence of
ordinals $\left(  \rho_{k}\right)  _{k=0}^{\infty}$ with $\rho_{0}=0$ that
strictly increases to $\omega^{\xi_{1}}.$ Let $\lambda$ be any ordinal that is
less than $\omega^{\xi_{2}}.$ Fix a function $u:[0,\omega^{\lambda
})\rightarrow\mathbb{N}$ such that $\left\{  \alpha\in\lbrack0,\omega
^{\lambda}):u\left(  \alpha\right)  \leq k\right\}  $ is finite for all
$k\in\mathbb{N}.$ Define real-valued functions $f$ and $g$ on $K$ as follows.
If $t\in K_{\omega^{\xi_{1}}\cdot\lambda},$ let $f\left(  t\right)  =g\left(
t\right)  =0.$ If $t\in$ $K_{\omega^{\xi_{1}}\cdot\alpha+\rho_{k-1}}\setminus$
$K_{\omega^{\xi_{1}}\cdot\alpha+\rho_{k}}$ for some $\alpha<\omega^{\lambda}$
and $k\in\mathbb{N},$ let $f\left(  t\right)  =\dfrac{h\left(  t\right)
}{ku\left(  \alpha\right)  }$ and $g\left(  t\right)  =ku\left(
\alpha\right)  .$ Notice that $fg=h\chi_{V},$ where $V=K\setminus
K^{\omega^{\xi_{1}}\cdot\lambda}\left(  h,1\right)  .$ It follows from Lemma
\ref{P4} that $K^{\eta}\left(  h,1\right)  \setminus K^{\omega^{\xi_{1}}%
\cdot\lambda}\left(  h,1\right)  \subseteq K^{\eta}\left(  fg,1\right)  $ for
all $\eta<\omega_{1}.$ Since $K^{\omega^{\xi}}\left(  h,1\right)
\neq\emptyset,$ and $h\in\frak{B}_{1}\left(  K\right)  ,$ $K^{\eta}\left(
h,1\right)  \setminus K^{\omega^{\xi_{1}}\cdot\lambda}\left(  h,1\right)
\neq\emptyset$ for all $\eta<\omega^{\xi_{1}}\cdot\lambda.$ Thus $K^{\eta
}\left(  fg,1\right)  \neq\emptyset$ for all $\eta<\omega^{\xi_{1}}%
\cdot\lambda.$ Hence $\beta\left(  fg\right)  \geq\omega^{\xi_{1}}\cdot\lambda.$

We now turn to the calculation of $\beta\left(  g\right)  $ and $\beta\left(
f\right)  .$ First notice that the sets $K_{\omega^{\xi_{1}}\cdot\alpha
+\rho_{k-1}}\setminus K_{\omega^{\xi_{1}}\cdot\alpha+\rho_{k}}$,
${k\in\mathbb{N}}$, form a partition of $K_{\omega^{\xi_{1}}\cdot\alpha
}\setminus$ $K_{\omega^{\xi_{1}}\cdot\left(  \alpha+1\right)  }$ into
relatively open sets for any $\alpha<\omega^{\lambda},$ and that $g$ is
constant on each set $K_{\omega^{\xi_{1}}\cdot\alpha+\rho_{k-1}}\setminus
K_{\omega^{\xi_{1}}\cdot\alpha+\rho_{k}}.$ Hence the restriction of $g$ to
$K_{\omega^{\xi_{1}}\cdot\alpha}\setminus$ $K_{\omega^{\xi_{1}}\cdot\left(
\alpha+1\right)  }$ is a continuous function for each $\alpha<\omega^{\lambda
}.$ It follows readily by induction that for any $\varepsilon>0,$
$K^{\alpha}\left(  g,\varepsilon\right)  \subseteq K_{\omega^{\xi_{1}}%
\cdot\alpha}$ for all $\alpha\leq\omega^{\lambda}.$ But $g=0$ on
$K_{\omega^{\xi_{1}}\cdot\alpha}.$ Thus $K^{\omega^{\lambda}+1}\left(
g,\varepsilon\right)  =\emptyset.$ Therefore $\beta\left(  g\right)
\leq\omega^{\lambda}+1\leq\omega^{\xi_{2}}.$

Finally, consider the function $f$. Let $k_{0}\in\mathbb{N}$ be given. The
set
\[
A=\left\{  \left(  \alpha,k\right)  :k\in\mathbb{N},\text{ }\alpha\in
\lbrack0,\omega^{\lambda}),\,ku\left(  \alpha\right)  \leq k_{0}\right\}
\]
is finite. List the elements of $A$ in a finite sequence $\left(  \left(
\alpha_{i},k_{i}\right)  \right)  _{i=1}^{j}$ in lexicographical order. Then
$\left|  f\left(  t_{1}\right)  -f\left(  t_{2}\right)  \right|  <\dfrac
{1}{k_{0}}$ for all $t_{1},$\thinspace$t_{2}\in K\setminus K_{\omega^{\xi_{1}%
}\cdot\alpha_{1}+\rho_{k_{1}-1}}.$ Hence $K^{1}\left(  f,\tfrac{1}{k_{0}%
}\right)  \subseteq K_{\omega^{\xi_{1}}\cdot\alpha_{1}+\rho_{k_{1}-1}}.$ Note
that $f=\frac{h}{k_{1}u\left(  \alpha_{1}\right)  }$ on $K_{\omega^{\xi_{1}%
}\cdot\alpha_{1}+\rho_{k_{1}-1}}\setminus$ $K_{\omega^{\xi_{1}}\cdot\alpha
_{1}+\rho_{k_{1}}}.$ Thus $K^{1+\eta}\left(  f,\tfrac{1}{k_{0}}\right)
\subseteq K_{\omega^{\xi_{1}}\cdot\alpha_{1}+\rho_{k_{1}-1}+\eta}$ for all
$\eta$ such that $\omega^{\xi_{1}}\cdot\alpha_{1}+\rho_{k_{1}-1}+\eta
\leq\omega^{\xi_{1}}\cdot\alpha_{1}+\rho_{k_{1}}.$ Let $\eta_{0}$ be such that
$\omega^{\xi_{1}}\cdot\alpha_{1}+\rho_{k_{1}-1}+\eta_{0}=\omega^{\xi_{1}}%
\cdot\alpha_{1}+\rho_{k_{1}}.$ Then $\eta_{0}\leq\rho_{k_{1}}.$ Therefore,%

\[
K^{1+\rho_{k_{1}}}\left(  f,\tfrac{1}{k_{0}}\right)  \subseteq K^{1+\eta_{0}%
}\left(  f,\tfrac{1}{k_{0}}\right)  \subseteq K_{\omega^{\xi_{1}}\cdot
\alpha_{1}+\rho_{k_{1}}}.
\]
Repeating the argument, we see that
\[
K^{\rho}\left(  f,\dfrac{1}{k_{j}}\right)  \subseteq K_{\omega^{\xi_{1}}%
\cdot\alpha+\rho_{k_{j}}},
\]
where $\rho=1+\rho_{k_{1}}+1+\rho_{k_{2}}+...+1+\rho_{k_{j}}.$ Since $0\leq
f\left(  t\right)  <\dfrac{1}{k_{0}}$ for all $t\in K_{\omega^{\xi_{1}}%
\cdot\alpha+\rho_{k_{j}}},$%
\[
K^{\rho+1}\left(  f,\dfrac{1}{k_{j}}\right)  =\emptyset.
\]
As $\left(  \rho_{k}\right)  $ increases to $\omega^{\xi_{1}},$ $\rho
+1<\omega^{\xi_{1}}.$ Hence $K^{\omega^{\xi_{1}}}\left(  f,\dfrac{1}{k_{0}%
}\right)  =\emptyset$ for any $k_{0}\in\mathbb{N}.$ It follows that
$\beta\left(  f\right)  \leq\omega^{\xi_{1}}.$ Summarizing, we have functions
$f$ and $g$ such that $f\in\frak{B}_{1}^{\xi_{1}}\left(  K\right)  $,
$g\in\frak{B}_{1}^{\xi_{2}}\left(  K\right)  $ and $\beta\left(  fg\right)
\geq\omega^{\xi_{1}}\cdot\lambda.$ Since $\lambda<\omega^{\xi_{2}}$ is
arbitrary, the theorem is proved.
\end{proof}


\begin{thebibliography}{99}
\bibitem{D}\textsc{James Dugundji}, Topology, Allyn and Bacon, Inc., Boston, 1966.

\bibitem {H-O-R}\textsc{R. Haydon, E. Odell and H. P. Rosenthal}, On certain
classes of Baire-1 functions with applications to Banach space theory,
Functional Analysis Proceedings, The University of Texas at Austin 1987-89,
Lecture Notes in Math., vol. 1470, Springer-Verlag, New York, 1991, pp. 1-35.

\bibitem {K-L}\textsc{A. S. Kechris and A. Louveau}, A classification of
Baire-1 functions, Trans. Amer. Math. Soc. \textbf{318}(1990), 209-236.

\bibitem {K}\textsc{P. Kiriakouli}, A classification of Baire-1 functions,
Trans. Amer. Math. Soc. \textbf{351}(1999), 4599-4609.
\end{thebibliography}
\end{document}